\documentclass
[12pt]{amsart}
\usepackage{amssymb}

\usepackage[all]{xy}

\newcommand{\Q}{{\mathbb {Q}}}

\newcommand{\G}{{\bf{G}}}
\newcommand{\T}{{\bf{T}}}
\newcommand{\U}{{\bf{U}}}
\newcommand{\V}{{\bf{V}}}
\newcommand{\R}{{\mathbb{R}}}
\newcommand{\Z}{{\mathbb{Z}}}
\newcommand{\C}{{\mathbb{C}}}
\newcommand{\h}{{\bf{H}}}

\newcommand{\se}{{\bf{S}}}

\newcommand{\te}{{\bf{T}}}

\newcommand{\ve}{{\bf{V}}}
\newcommand{\pe}{{\bf{P}}}

\newcommand{\N}{{\mathbb{N}}}

\newcommand{\SL}{\operatorname{SL}}

\newcommand{\Lie}{\operatorname{Lie}}

\newcommand{\OO}{{\mathcal O}}

\newcommand{\LL}{{\mbox{\boldmath $\mathfrak{g}$}}}

\newcommand{\rank}{{\rm rank}}

\def\bydefn{\stackrel{def}{=}}

\theoremstyle{plain}
\newtheorem{thm}{Theorem}[section]
\newtheorem{lem}[thm]{Lemma}
\newtheorem{prop}[thm]{Proposition}
\newtheorem{cor}[thm]{Corollary}

\theoremstyle{definition}

\def\bydefn{\stackrel{def}{=}}

\begin{document}

\title[Locally divergent orbits]{Closures of locally divergent orbits of maximal tori and values of homogeneous forms}

\author{George Tomanov}
\address{Institut Camille Jordan, Universit\'e Claude Bernard - Lyon
I, B\^atiment de Math\'ematiques, 43, Bld. du 11 Novembre 1918,
69622 Villeurbanne Cedex, France {\tt tomanov@math.univ-lyon1.fr}}


\maketitle

\begin{abstract} Let $\G$ be a semisimple algebraic group over a number field $K$, $\mathcal{S}$ a finite set of places of $K$,
$K_\mathcal{S}$ the direct product of the completions $K_v, v \in \mathcal{S}$, and $\OO$ the ring of $\mathcal{S}$-integers of $K$.
Let $G = \G(K_\mathcal{S})$, $\Gamma = \G(\OO)$ and $\pi:G \rightarrow G/\Gamma$ the quotient map.  We describe the closures of the locally divergent
orbits ${T\pi(g)}$ 
where $T$ is a maximal $K_\mathcal{S}$-split torus
in $G$. If $\# S = 2$ then the closure
$\overline{T\pi(g)}$ is a finite union of $T$-orbits stratified in terms of parabolic subgroups of $\G \times \G$ and, consequently, $\overline{T\pi(g)}$
is homogeneous (i.e., $\overline{T\pi(g)}= H\pi(g)$ for a subgroup $H$ of $G$) if and only if ${T\pi(g)}$ is closed.
On the other hand, if $\# \mathcal{S} > 2$ and $K$ is not a $\mathrm{CM}$-field then $\overline{T\pi(g)}$ is homogeneous for $\G = \mathbf{SL}_{n}$
and, generally, non-homogeneous but squeezed between closed orbits of
two reductive subgroups of equal semisimple $K$-ranks for $\G \neq \mathbf{SL}_{n}$.
As an application, we prove that $\overline{f(\OO^{n})} = K_{\mathcal{S}}$ for the class of non-rational locally $K$-decomposable homogeneous
forms $f \in K_{\mathcal{S}}[x_1, \cdots, x_{n}]$.
\end{abstract}

\section{Introduction} \label{Introduction}
Let $\G$ be a semisimple algebraic group defined over a number field $K$. Let $\mathcal{S}$ be a finite set of
places of $K$
containing the archimedean ones and let $\OO$ be the ring of $\mathcal{S}$-integers in $K$. Denote by $K_v$, $v \in \mathcal{S}$,
the completion of $K$ with respect to $v$ and by $K_{\mathcal{S}}$ the direct product of the topological fields $K_v$. 
 Put $G = \G(K_{\mathcal{S}})$.
The group  $G$ is naturally identified with the direct product of the locally compact groups $G_v = \G(K_v)$, $v \in \mathcal{S}$,
and $\G(K)$ is diagonally imbedded in $G$.
 Let $\Gamma$ be an $\mathcal{S}$-arithmetic subgroup of $G$, that is, $\Gamma \cap \G(\OO)$ have
 finite index in both $\Gamma$ and $\G(\OO)$. Recall that the homogeneous space $G/\Gamma$ endowed with the quotient topology has
 finite volume with respect to the Haar measure.
 Let $H$ be a closed subgroup of $G$ acting on
 $G/\Gamma$ by left translations, that is, $$
h\pi(g) \bydefn \pi(hg), \forall h \in H,
$$
where $\pi: G \rightarrow G/\Gamma$ is the quotient map.
An orbit
$H\pi(g)$ is called {\it divergent} if the orbit map $H \rightarrow G/\Gamma, h
\mapsto h\pi(g)$, is proper, i.e., if $\{h_i\pi(g)\}$ leaves
compacts of $G/\Gamma$ whenever $\{h_i\}$ leaves compacts of $H$. It is clear that the divergent orbits are closed.
The closure $\overline{H\pi(g)}$ of ${H\pi(g)}$ in $G/\Gamma$ is called \textit{homogeneous} if $\overline{H\pi(g)} = L\pi(g)$ for
a closed subgroup $L$ of $G$.

Fix a maximal $K$-split
torus $\te$ of $\G$ and for every $v \in \mathcal{S}$ a maximal $K_v$-split torus $\te_v$ of $\G$ containing $\te$.
Recall that, given a field extension $F/K$, the $F$-rank of $\G$, denoted by $\mathrm{rank}_F \G$, is the common dimension
of the maximal $F$-split tori of $\G$.
So, $\rank_{K_v} \G \geq \rank_K \G$ and $\rank_{K_v} \G = \rank_K \G$ if and only if $\te = \te_v$. Let $T_v = \te_v(K_v)$ and
$T = \prod_{v \in \mathcal{S}}T_v \subset G$. 
 An orbit $T \pi(g)$
is called $\it{locally \  divergent}$ if $T_{v}\pi(g)$ is divergent for every $v \in \mathcal{S}$.

The locally divergent orbits, in general, and the closed locally divergent orbits, in particular, are completely described by the following.

\begin{thm} \mbox{\rm{(}}\cite[Theorem 1.4 \ and Corollary 1.5]{Toma}\mbox{\rm{)}}
\label{ldo}
With the above notation, we have:
\begin{enumerate}
\item[(a)] An orbit $T_v\pi(g)$ is divergent if and
only if
\begin{equation}
\label{rank}
\rank_{K_v} \G = \rank_K \G
\end{equation}
and
\begin{equation}
\label{rank'}
g \in \mathcal{N}_{G}(T_{v})\G(K),
\end{equation}
where $\mathcal{N}_{G}(T_v)$ is  the normalizer of $T_v$ in $G$. So, $T\pi(g)$ is locally divergent if and
only if (\ref{rank}) and (\ref{rank'}) hold for all $v \in \mathcal{S}$;
\item[(b)]  An orbit $T\pi(g)$ is both locally divergent orbit  and closed if and
only if (\ref{rank}) holds for all $v \in \mathcal{S}$ and
\begin{equation*}
g \in \mathcal{N}_{G}(T)\G(K),
\end{equation*}
where $\mathcal{N}_{G}(T)$ is  the normalizer of $T$ in $G$.
\end{enumerate}
\end{thm}

Our Theorem \ref{ldo} is the accomplishment of several works (cf. \cite{Tomanov-Weiss}, \cite{Weiss1} and  \cite{Toma2}),
the first being the classification by G.A.Margulis (see \cite[Appendix]{Tomanov-Weiss}) of
the divergent orbits for the action of the group of diagonal matrices in $\mathrm{SL}_n(\R)$
on $\mathrm{SL}_n(\R)/\mathrm{SL}_n(\Z)$.
In the present paper, using completely different ideas,
 we describe the closures of the locally divergent non-closed
 orbits. According to Theorem \ref{ldo} such orbits exist if $\rank_K \G > 0$, $\#\mathcal{S} \geq 2$ and (\ref{rank}) holds for all $v \in \mathcal{S}$.

Essentially due to applications
to Diophantine approximation of numbers, the study of orbit closures in $G/\Gamma$ for different kind of subgroups $H$ of $G$ attracted considerable
interest during the last decades. In view of the classical result \cite{Margulis non-divergence},
the orbits of the $1$-parameter unipotent subgroups are always recurrent. Hence
 if $H$ is generated by $1$-parameter unipotent subgroups then $H\pi(g)$ is never divergent.
Moreover, for such kind of $H$
it is proved by M.Ratner in \cite{Ratner} and \cite{Ratner+},
in the real setting, and in
\cite{Mar-Tom} and  \cite{Ratner++} (see also \cite{Mar-Tom+} and \cite{Toma4}), in the $\mathcal{S}$-adic setting, that $\overline{H\pi(g)}$
is homogeneous.  The special case when $H = \mathrm{SO}(q)$, where $q$ is a non-degenerate indefinite quadratic form on
$\R^n, n \geq 3$, is acting on $\mathrm{SL}_n(\R)/\mathrm{SL}_n(\Z)$ was first established by Margulis \cite[Theorem 2]{Margulis Oppenheim 2}
for bounded orbits and by Dani and Margulis \cite[Theorem 2]{Dani-Margulis1}, in general. The latter  implies
that $\overline{q(\Z^n)} = \R$ provided $q$ is not a multiple of a form with integer coefficients (see \cite[Theorem 1]{Dani-Margulis1})
strengthening \cite[Theorem 1']{Margulis Oppenheim 2} which confirms the A.Oppenheim conjecture. Also by using the homogeneous space approach,
the $\mathcal{S}$-adic version of the Oppenheim conjecture is proved by A.Borel and G.Prasad \cite{Borel-Prasad}.
The dynamics of the action of split tori $T \subset G$ on $G/\Gamma$ is much less understood and reveals completely different phenomena.
Concerning the orbit closures, it was believed up to recently
 that $\overline{T\pi(g)}$ is homogeneous if $G/\Gamma$ does not admit rank 1 $T$-invariant
factors (see \cite[Conjecture 1]{Margulis problems survey}). Affirmative results in the simplest case when $G = \mathrm{SL}_2(K_1) \times \mathrm{SL}_2(K_2)$,
where $K_1$ and $K_2$ are local fields, $\Gamma$ is an irreducible lattice in $G$ and $T$ is the direct product of the subgroups of
diagonal matrices in the first  and the second copy of $\mathrm{SL}_2$ have been obtained in \cite{F} and \cite{Mozes}).
Nevertheless, it turned out that if $T\pi(g)$ is locally divergent
then $\overline{T\pi(g)}$ is
homogeneous only if $T\pi(g)$ is closed which, in view of Theorem \ref{ldo}(b), contradicts \cite[Conjecture 1]{Margulis problems survey}
(cf. \cite[Corollary 1.2]{Toma3}).
The result is generalized and strengthened for
 arbitrary semisimple groups by Theorems \ref{r=2} and \ref{THM2} of the present paper.
Sparse examples of non-homogeneous orbit closures of completely different nature
are given in \cite{Mau} for the action of a $n-2$-dimensional split torus on $\mathrm{SL}_n(\R)/\Gamma$, $n \geq 6$,
and in \cite{Shapira} and \cite{L-Sha} for the action of a $2$-dimensional split torus on $\mathrm{SL}_3(\R)/\mathrm{SL}_3(\Z)$.
The understanding of orbit closures of maximal split tori admits deep number theoretical applications.
For instance, if $f \in \R[x_1, \cdots, x_{n}]$ is a product
of $n \geq 3$ linearly independent real linear forms then \cite[Conjecture 8]{Margulis problems survey} claims that $f$ is
a multiple of a form with integer coefficients whenever $0$ is an isolated point in ${\overline{f(\Z^{n})}}$ and $f(\vec{a}) \neq 0$
for all $\vec{a} \in \Z^{n} \setminus \{\vec{0}\}$.
In terms of group actions, \cite[Conjecture 8]{Margulis problems survey} is equivalent to \cite[Conjecture 9]{Margulis problems survey} stating
that every bounded orbit for the action of the group of diagonal matrices on $\mathrm{SL}_n(\R)/\mathrm{SL}_n(\Z)$, $n \geq 3,$ is homogeneous.
From its side, \cite[Conjecture 9]{Margulis problems survey} implies a well-known conjecture of Littlewood \cite[Conjecture 7]{Margulis problems survey}
formulated around 1930 and seemingly still far from its final solution. (See \cite{EiKaLi}, \cite{EiKl} and \cite{EinLind} for
recent results on the Littlewood conjecture.)
Along the same line, Theorem \ref{cor5} below implies that $\overline{f(\OO^n)} = K_\mathcal{S}$
for a natural class of non-rational forms $f$ on $K_\mathcal{S}^n$ (Theorem \ref{application}).

From now on, with the notation of Theorem \ref{ldo}, we suppose that $\#\mathcal{S} \geq 2$ and ${T\pi(g)}$ is a locally divergent orbit.
The cases $\#\mathcal{S} = 2$ and $\#\mathcal{S} > 2$ behave in drastically different ways.
The next two theorems describe $\overline{{T\pi(g)}}$ in both cases.

\begin{thm}
\label{r=2} Let $\#\mathcal{S} = 2$. Then
\begin{enumerate}
\item ${\overline{T\pi(g)}}$ is a union  of finitely many $T$-orbits which are all locally divergent  and stratified in terms of parabolic subgroups of
$\G \times \G$;
\item ${T\pi(g)}$ is open in ${\overline{T\pi(g)}}$;
\item The following conditions are equivalent:
\begin{enumerate}
\item[(a)] ${T\pi(g)}$ is closed,
\item[(b)] $\overline{T\pi(g)}$ is homogenous,
\item[(c)] $g \in \mathcal{N}_G(T)\G(K)$.
\end{enumerate}
\end{enumerate}
\end{thm}

Theorem \ref{r=2} is a particular case of stronger but more
technically formulated results proven in \S5. More precisely,
the part (1) of Theorem \ref{r=2} is a particular case of Theorem \ref{THM2}, its part (2) is a particular case of Corollary \ref{cor1} and its part (3)
coincides with Corollary \ref{cor4}. The $T$-orbits contained in ${\overline{T\pi(g)}}$ are stratified in the following sense. (See \S5 for details.)
Given a locally divergent orbit $T\pi(g)$, we define
 a finite set $\mathcal{P}(g)$ of parabolic subgroups of $\G \times \G$ and associate to each $\pe \in \mathcal{P}(g)$  a $T$-orbit
 $\mathrm{Orb}_{g}(\pe)$ contained in ${\overline{T\pi(g)}}$. We have $\G \times \G \in \mathcal{P}(g)$ and $T\pi(g) = \mathrm{Orb}_{g}(\G \times \G)$.
If $\pe \in \mathcal{P}(g)$ then $\overline{\mathrm{Orb}_{g}(\pe)} = \bigcup_{\pe' \in \mathcal{P}(g), \ \pe' \subset \pe} \mathrm{Orb}_{g}(\pe')$, in
particular, ${\overline{T\pi(g)}} = \bigcup_{\pe \in \mathcal{P}(g)} \mathrm{Orb}_{g}(\pe)$ and the closed
 $T$-orbits in ${\overline{T\pi(g)}}$ correspond to the minimal parabolic subgroups contained in $\mathcal{P}(g)$ (see Corollary \ref{cor2}(b)).
 The correspondence between the parabolic
 subgroups and the $T$-orbits in ${\overline{T\pi(g)}}$ becomes bijective under Zariski topology density conditions on $g \in G$ (see Corollary \ref{cor3}).

Recall that
the semi-simple $K$-rank of a reductive $K$-group $\mathbf{H}$, denoted by $\mathrm{s.s.rank}_K (\mathbf{H})$, is equal to
$\rank_K \mathcal{D}(\mathbf{H})$ where $\mathcal{D}(\mathbf{H})$ is  the derived subgroup of $\mathbf{H}$. Also, $K$ is called a $\mathrm{CM}$-field
if it is a quadratic extension $K/F$ where $F$ is a totally real
number field but $K$ is totally imaginary. So, a totally real
number field is not a $\mathrm{CM}$-field.

The main result for $\#\mathcal{S} > 2$ is the following.
\begin{thm}
\label{r>2}
Let  $\#\mathcal{S} > 2$ and  $K$ be not a $\mathrm{CM}$-field.
Then there exist
$h_1$ and $h_2 \in \mathcal{N}_G(T)\G(K)$ and reductive $K$-subgroups $\mathbf{H}_1$ and $\mathbf{H}_2$ of $\G$ such that
$\mathbf{H}_1 \subset \mathbf{H}_2$, $\mathrm{rank}_K (\mathbf{H}_1) = \mathrm{rank}_K (\mathbf{H}_2) = \mathrm{rank}_K (\mathbf{G})$,
\begin{equation}
\label{rank+}
\mathrm{s.s.rank}_K (\mathbf{H}_1) =\mathrm{s.s.rank}_K (\mathbf{H}_2),
\end{equation}
and
\begin{equation}
\label{skeezing}
h_2H_2 \pi(e) \subseteq \overline{T\pi(g)} \subseteq h_1H_1\pi(e),
\end{equation}
where $H_1 = \mathbf{H}_1(K_{\mathcal{S}})$, $H_2$ is a subgroup of finite index in $\mathbf{H}_2(K_{\mathcal{S}})$, and the orbits
$h_1H_1\pi(e)$ and $h_2H_2 \pi(e)$ are closed and $T$-invariant.
\end{thm}

In the important case
$\G = \mathbf{SL}_{n}$, Theorem \ref{r>2} implies

\begin{thm}
\label{cor5} Let $\G = \mathbf{SL}_{n}$, $\#\mathcal{S} > 2$ and  $K$ be not a $\mathrm{CM}$-field.
Then $\overline{T\pi(g)} = H\pi(g)$,
where $H$ is a closed subgroup of $G$.
\end{thm}

Theorem \ref{contre-example}, proven in \S8, provides examples showing that
Theorem \ref{cor5} is not valid for $\mathrm{CM}$-fields and, also,
that $\overline{T\pi(g)}$ as in the formulation of Theorem \ref{r>2} might not be homogeneous.
The orbit ${T\pi(g)}$ given by Theorem \ref{contre-example} is such that
$\overline{T\pi(g)} \setminus {T\pi(g)}$ is not contained in a countable union of closed orbits of proper subgroups of $G$ in contrast to
the orbits ${T\pi(g)}$ with non-homogeneous closures given in \cite{Mau}, \cite{Shapira}, \cite{L-Sha} and our Theorem \ref{r=2}
where $\overline{T\pi(g)} \setminus {T\pi(g)}$ is always contained in a finite union of closed orbits of proper subgroups of $G$.

Before stating the number theoretical application of  Theorem \ref{cor5} we need to set up some notation and formulate a general conjecture.
As usual,  $K_{\mathcal{S}}[\,\vec{x}\,]$ denotes the ring of polynomials in $n$
variables $\vec{x} = (x_1,\ldots,x_{n})$ with coefficients from the ring $K_{\mathcal{S}}$. We suppose that $n \geq 2$.
Note that $K_{\mathcal{S}}[ \vec{x}] = \underset{v \in \mathcal{S}}{\prod} K_v[ \vec{x}]$ and the ring $K[ \vec{x}]$
is identified with its diagonal imbedding in $K_{\mathcal{S}}[ \vec{x}]$.
Further on, $f(\vec{x}) = (f_v(\vec{x}))_{v \in \mathcal{S}} \in K_{\mathcal{S}}[ \vec{x}]$ is
a decomposable  (over $K_{\mathcal{S}}$) non-degenerated homogeneous form, that is,
$f(\vec{x}) = l_{1}(\vec{x}) \cdots l_{m}(\vec{x})$, where $l_{1}(\vec{x}), \ldots,
l_{m}(\vec{x})$ are linearly independent over $K_{\mathcal{S}}$ linear forms with coefficients from $K_{\mathcal{S}}$.  Equivalently,
we suppose that every $f_v(\vec{x}) = l_{1}^{(v)}(\vec{x}) \cdots l_{m}^{(v)}(\vec{x})$, where $l_{1}^{(v)}(\vec{x}), \ldots,
l_{m}^{(v)}(\vec{x})$ are linearly independent over $K_v$ linear forms with coefficients from $K_{v}$.
The form $f$ is called rational if $f(\vec{x}) = c \cdot h(\vec{x})$, where
$h(\vec{x}) \in K[\vec{x}]$ and $c \in K_{\mathcal{S}}$, and non-rational, otherwise. According to \cite[Theorem 1.8]{Toma} $f$ is rational
if and only if $f(\mathcal{O}^n)$ is discrete in $K_{\mathcal{S}}$. For non-rational forms $f$
the following conjecture is plausible.

\medskip

\textbf{Conjecture 1.} Suppose that $\#\mathcal{S} > 2$, $K$ is not a $\mathrm{CM}$-field and $f$ is non-rational.
Then $\overline{f(\OO^n)} = K_{\mathcal{S}}$.

\medskip

The form $f$ is called locally $K$-decomposable if for every $v \in \mathcal{S}$ each of the linear forms $l_{1}^{(v)}(\vec{x}), \ldots,
l_{m}^{(v)}(\vec{x})$ is a multiple of a linear form with coefficients from $K$. Theorem \ref{cor5} implies:

\begin{thm}
\label{application} Conjecture 1 is true for the locally $K$-decomposable homogeneous forms.
\end{thm}

Theorems \ref{two valuations} and \ref{example CM field} of section 8.3 show that the analog of Theorem \ref{application} (and, therefore, of
Conjecture 1)
is not true if $\#\mathcal{S} = 2$ or $\#\mathcal{S} \geq 2$ and $K$ is a $\mathrm{CM}$-field.

\section{Preliminaries: notation and some basic concepts} \label{section:prelims}

\subsection{Numbers}

As usual, $\N$, $\Z$, $\Q$, $\R$, and $\C$ denote the non-negative integer, integer, rational, real and complex
numbers, respectively. Also, $\N_+ = \{x \in \N : x > 0\}$ and  $\R_+ = \{x \in \R : x > 0\}$. 

In this paper ${K}$ is a number field, that is, a finite
extension of $\Q$.
If $v$ is a place of $K$
then ${K}_v$ is the completion of $K$ with respect to $v$ and $| \cdot|_v$ is the corresponding \textit{normalized} norm on $K_v$
(see \cite[ch.2, \S7]{Cassels-Fröhlich}).
Recall that if $K_v = \R$ (respectively, $K_v = \C$) then $| \cdot|_v$ is the absolute value on $\R$ (respectively, the square of the absolute value on $\C$).
If $v$ is non-archimedean then ${\OO}_v = \{x \in K_v : |\
 x\ |_v \leq 1 \}$ denotes the ring of integers in ${K}_v$.

Further on we denote by $\hat{K}$ an universal domain, that is, an algebraically closed
field containing $K$ and all completions of $K_v$ of $K$.

We fix a finite set $\mathcal{S} = \{v_1, \cdots, v_r\}$ of places of ${K}$ containing all
archimedean places of $K$. The archimedean places in $\mathcal{S}$ will be denoted by
$\mathcal{S}_\infty$. We let $\mathcal{S}_f = \mathcal{S} \setminus \mathcal{S}_\infty$.

Sometimes we will write $K_i$ instead of $K_{v_i}$ and $|\cdot|_i$ instead of $|\cdot|_{v_i}$.

We denote by $\mathcal{O}$ the ring of $\mathcal{S}$-integers in $K$, i.e.,
$\mathcal{O} = K \bigcap (\bigcap_{v \notin \mathcal{S}}\mathcal{O}_v)$. Also,
$\mathcal{O}_{\infty} = K \bigcap (\bigcap_{v \notin \mathcal{S}_\infty}\mathcal{O}_v)$ is the ring of
integers in $K$.

Let ${K}_\mathcal{S} \bydefn \prod_{v \in
\mathcal{S}}{K}_v$ considered with the product topology. The diagonal imbedding of  $K$ into the topological ring ${K}_\mathcal{S}$ is dense
 and $\mathcal{O}$ is a lattice in $K_{\mathcal{S}}$. We denote $K_\infty = \prod_{v \in
\mathcal{S}_\infty}{K}_v$.

As usual, if $R$ is a ring $R^*$ denotes
the multiplicative group of units of $R$.

\subsection{Groups.}
Further on, we use boldface letters to denote the algebraic groups defined over $K$ (shortly, the $K$-algebraic groups or the algebraic $K$-groups).
Let $\h$ be a $K$-algebraic group.
As usual, $\mathcal{R}_u(\h)$ (respectively, $\mathrm{Lie}(\G)$) stands for the unipotent radical (respectively, the Lie algebra) of $\h$.
Given $v \in \mathcal{S}$, we write $H_v \bydefn \h(K_v)$ or simply $H_i$ if $\mathcal{S} = \{v_1, \cdots, v_r\}$ and $v = v_i$.
Put $H \bydefn \h(K_{\mathcal{S}})$. The group $\h(K)$ is identified with its diagonal imbedding in $H$. On every $H_v$ we have \textit{Zariski topology} induced by the Zariski topology on $\h$ and
\textit{Hausdorff topology} induced by the Hausdorff topology on $K_v$. The formal product
of the Zariski (resp., the Hausdorff) topologies on $H_v$, $v
\in S$, is the Zariski (respectively, the Hausdorff) topology on
$H$. In order to distinguish the two topologies, all topological
notions connected with the first one will be used with the prefix
"Zariski".

The algebraic groups in this paper are always linear. Every $K$-algebraic group $\h$ is a Zariski closed $K$-subgroup of
$\mathbf{GL}_l$ for some $l \in \N_+$. The group $\mathbf{GL}_l$ itself is identified with
$\mathrm{GL}_l(\hat{K})$ where $\hat{K}$ is the universal domain defined in 2.1. We have $\mathbf{GL}_1(\OO) = \OO^*$ and
$\h(\OO) = \mathbf{GL}_l(\OO) \cap \h$. A $K$-subgroup $\te$ of $\h$ is a $K$-split torus if $\te$ is $K$-isomorphic to
$\mathbf{GL}_1^{d}$ for some $d \in \mathbb{N}$.
A subgroup $\Delta$ of $\h$ is called $\mathcal{S}$-\textit{arithmetic} if
 $\Delta$ and $\h(\OO)$ are commensurable, that is, if
$\Delta \cap \h(\OO)$ has finite index in both $\Delta$ and
$\h(\OO)$. Recall that if $\h$ is semisimple then $\Delta$ is a lattice in $H$, i.e. $H/\Delta$ has finite Haar measure.

The Zariski connected component of the identity $e \in \mathbf{H}$ is denoted by $\mathbf{H}^\circ$.
In the case of a real Lie group $L$ the connected component of the identity is denoted by $L^\bullet $.

If $A$ and $B$ are subgroups of an abstract group $G$ then $\mathcal{N}_{A}{({B})}$
(resp., $\mathcal{Z}_A({B})$) is the normalizer
(resp., the centralizer) of ${B}$ in $A$. As usual, $\mathcal{Z}(G)$ denotes the center of $G$ and $\mathcal{D}(G)$ the derived subgroup of $G$.

\subsection{$K$-roots.}
In this paper $\G$ is a connected, semisimple, $K$-isotropic algebraic group and $\te$ is a maximal $K$-split
torus in $\G$.

We denote by   $\Phi (\equiv \Phi(\te,\G))$ the system of $K$-roots with respect to  $\mathbf{T}$.
Let $\Phi^{+}$ be a system of positive $K$-roots in $\Phi$
 and $\Pi$ be the system of simple roots in $\Phi^{+}$.
(We refer to
 \cite[\S 21.1]{Borel} for the standard definitions related to the $K$-roots.)

If $\chi \in \Phi$ we let $\LL_{\chi}$ be the corresponding root-space in  $\mathrm{Lie}(\G)$.
For every $\alpha \in \Pi$ we define a projection $\pi_\alpha: \Phi \rightarrow \Z$ by $\pi_\alpha(\chi) = n_\alpha$
where $\chi = \sum_{\beta \in \Pi}n_\beta \beta$.

Let $\Psi \subset \Pi$ and $\mathbf{T}_{\Psi} \bydefn (\bigcap_{\alpha \in \Psi} \ker(\alpha))^\circ$.
We denote by $\mathbf{P_\Psi}$ the (standard) parabolic subgroup corresponding to $\Psi$
and by $\mathbf{P^-_\Psi}$ the opposite parabolic subgroup corresponding to $\Psi$. The centralizer $\mathcal{Z}_\G(\mathbf{T}_{\Psi})$
is a common Levi subgroup of $\mathbf{P_\Psi}$ and $\mathbf{P^-_\Psi}$, $\mathbf{P_\Psi} = \mathcal{Z}_\G(\mathbf{T}_{\Psi}) \ltimes \mathcal{R}_u(\mathbf{P_\Psi})$  and
$\mathbf{P^-_\Psi} = \mathcal{Z}_\G(\mathbf{T}_{\Psi}) \ltimes \mathcal{R}_u(\mathbf{P^-_\Psi})$. We will often use the
simpler notation $\mathbf{V_\Psi} \bydefn \mathcal{R}_u(\mathbf{P}_\Psi)$ and $\mathbf{V^-_\Psi} \bydefn \mathcal{R}_u(\mathbf{P}^-_\Psi)$.
Recall that
\begin{equation}
\label{eq: uniprad of max parabolic}
\Lie (\mathbf{V_\Psi}) =  \bigoplus_{\exists \alpha \in
\Pi \setminus \Psi, \ \pi_{\alpha}(\chi) > 0} \LL_{\chi},
\end{equation}
\begin{equation}
\label{eq: uniprad- of max parabolic}
\Lie (\mathbf{V^-_\Psi}) =  \bigoplus_{\exists \alpha \in
\Pi \setminus \Psi, \ \pi_{\alpha}(\chi) < 0} \LL_{\chi},
\end{equation}
and
\begin{equation}
\label{eq: Levi factor}
\Lie(\mathcal{Z}_\G(\mathbf{T}_\Psi)) = \Lie(\mathcal{Z}_\G(\mathbf{T})) \oplus
\bigoplus_{\forall \alpha \in \Pi \setminus \Psi, \ \pi_{\alpha}(\chi)=0}
\LL_{\chi}.
\end{equation}

It is well known that the map $\Psi \mapsto \mathbf{P}_\Psi$ is a bijection between the subsets of
$\Pi$ and the parabolic subgroups of $\mathbf{G}$
containing $\mathbf{B}$, cf. \cite[\S 21.11]{Borel}. Note that $\mathbf{P_\emptyset}$,
$\mathbf{P_\emptyset}^-$ are minimal parabolic subgroups and $\mathbf{G} = \mathbf{P}_\Pi = \mathbf{P}_\Pi^-$.

Given $\alpha \in \Phi$ we let $(\alpha)$ be the set of roots which are positive multiple of $\alpha$.
Then $\LL_{(\alpha)} \bydefn \bigoplus_{\beta \in (\alpha)} \LL_{\beta}$ is the Lie algebra of a unipotent group denoted by $\U_{(\alpha)}$.
Given $\Psi \subset \Pi$,  let $\Psi'$ be the set of all non-divisible positive roots $\chi$ such that $\exists \alpha \in
\Delta \setminus \Psi, \ \pi_{\alpha}(\chi) > 0$. Then the product morphism in any order
$\underset{\chi \in \Psi'}{\Pi}\U_{(\chi)} \rightarrow \ve_\Psi$ is an isomorphism of $K$-varieties, that is,  $\ve_\Psi$ is directly spanned
in any order by its subgroups $\U_{(\chi)}, \ \chi \in \Psi'$ \cite[21.9]{Borel}.

It follows from the above
definitions that $\Psi_1 \subset \Psi_2 \ \Leftrightarrow \ \mathbf{P}_{\Psi_1} \subset \mathbf{P}_{\Psi_2}
 \ \Leftrightarrow \ {\V}_{\Psi_1} \supset {\V}_{\Psi_2} \ \Leftrightarrow \ \mathcal{Z}_{\G}(\te_{\Psi_1})
  \subset \mathcal{Z}_{\G}(\te_{\Psi_2})$.
  Let $\V_{[\Psi_2 \setminus \Psi_1]} \bydefn
  \mathcal{Z}_\G({\T}_{\Psi_2}) \cap \V_{\Psi_1}$ and $\V^{-}_{[\Psi_2 \setminus \Psi_1]} \bydefn
  \mathcal{Z}_\G({\T}_{\Psi_2}) \cap \V^{-}_{\Psi_1}$.
It is easy to see that
\begin{equation}
\label{Bruhat}
\V_{\Psi_1} = \V_{\Psi_2}\V_{[\Psi_2 \setminus \Psi_1]} = \V_{[\Psi_2 \setminus \Psi_1]}\V_{\Psi_2}.
\end{equation}

Recall that the Weyl group $\mathcal{W} \bydefn \mathcal{N}_\G(\mathbf{T})/\mathcal{Z}_\G(\mathbf{T})$ acts by conjugation simply
transitively on the set of all minimal parabolic $K$-subgroups of $\G$ containing $\T$. When this does not lead to confusion,
we will identify the elements from $\mathcal{W}$ with their representatives from $\mathcal{N}_{\G}(\te)$. It is easy to see that
$\mathcal{W}_\Psi = \mathcal{N}_{\mathcal{Z}_\G(\mathbf{T}_\Psi)}(\mathbf{T})/\mathcal{Z}_\G(\mathbf{T})$ is the Weyl group of
$\mathcal{Z}_\G(\mathbf{T}_\Psi)$. Note that $\mathcal{W} = \mathcal{W}_{\emptyset}$.

We will denote by $\omega_0$ the element from $\mathcal{W}$ such that $\omega_0 \pe_{\emptyset} \omega_0^{-1} = \pe_{\emptyset}^{-}$.

\section{On the group of units of $\OO$}

 Recall that $\mathcal{S} = \{v_1, \cdots, v_r\}$, $r \geq 2$, $K_i = K_{v_i}$ and $K_{\mathcal{S}} = \underset{i}{\prod}K_i$.
By the $\mathcal{S}$-adic version of Dirichlet's unit theorem, the $\Z$-rank of $\OO^*$ is equal to $r - 1$.
Moreover, if $K_{\mathcal{S}}^1 = \{(x_1, \cdots, x_r) \in K_{\mathcal{S}}^*: |x_1|_1 \cdots |x_r|_r =
1\}$ then $\OO^*$ is a lattice of $K_{\mathcal{S}}^1$.

For every $m \in \mathbb{N}_+$, we denote $\OO^*_m = \{\xi^m | \xi
\in \OO^*\}$. The next proposition follows easily from the compactness of $K_{\mathcal{S}}^1/\OO^*_m$.

\begin{prop}
\label{integers1'} For every $m \in \mathbb{N}_+$ there exists a constant $\kappa_m > 1$ such that given $({a}_{i}) \in K_{\mathcal{S}}^1$
there exists $\xi \in \OO^*_m$ satisfying
$$
\frac{1}{\kappa_m} \leq |\xi {a}_{i}|_i \leq \kappa_m
$$
for all $1 \leq i \leq r$.
\end{prop}

Let $\mathcal{S}_\infty = \{v_1, \cdots, v_{r'}\}$ and
$\mathcal{S}_f = \{v_{r'+1}, \cdots, v_{r}\}$. So, $K_1 = \R$ or $\C$. Let $p: K_{\mathcal{S}}^1 \rightarrow K_1^*$ be the natural
projection and $L \bydefn \overline{p(\OO^*)}$. Remark that, in view of Dirichlet's unit theorem mentioned above, $L$ is co-compact in $K_1^*$.
Specific information on the
connected component $L^{\bullet}$ of $L$ is provided by the following

\begin{prop}
\label{integers1} With the above notation, we have:
\begin{enumerate}
\item If $r = 2$ then $L^{\bullet} = \{1\}$;
\item Let $r \geq 3$. We have:
\begin{enumerate}
\item $L^\bullet  \neq \{1\}$. In particular, if $K_1 = \R$ then $L^\bullet  = \R_+$;
\item Let $K_1 = \C$.
\begin{enumerate}
\item If $L^\bullet  = \R_+$ then $K$ is a $\mathrm{CM}$-field;
\item If $K$ is not a $\mathrm{CM}$-field and $L \neq \C^*$ then $L^\bullet $ coincides with the unit circle group in $\C^*$ unless $r = 3$
when $L^\bullet $ might be a spiral.
\end{enumerate}
\end{enumerate}
\end{enumerate}
\end{prop}
{\bf Proof.} $(1)$ follows easily from the compactness of $K_{\mathcal{S}}^1/\OO^*_m$.

$(2)$ If $r \geq 3$ in view of Dirichlet's unit theorem $\OO^*$ contains a subgroup of $\Z$-rank $2$. Therefore $p(\OO^*)$ is not discrete, proving that $L^\bullet  \neq \{1\}$. 

Let $K_1 = \C$. Suppose that $L^\bullet  = \R_+$. Therefore $L$ is a finite extension of $L^\bullet $. Hence there exists $m$ such that $p(\OO^*_m)$
is a dense subgroup of $L^\bullet $. Let $F$ be the number field generated over $\Q$ by $\OO^*_m$. Then $F$ is proper subfield of $K$ and its unit group
has the same $\Z$-rank as that of $K$, i.e. $K$ has a "unit defect". It is known that the fields with "unit defect" are exactly the $\mathrm{CM}$-fields  (cf.\cite{Remak}).

 It remains to consider the case when $K$ is not a $\mathrm{CM}$-field, $L \neq \C^*$ and $r > 3$. Since $L^\bullet $ is a $1$- dimensional subgroup of $\C^*$ we need to prove that $L^\bullet $ couldn't be a spiral.
 This will be deduced from the following six exponentials theorem due to Siegel: if $x_1, x_2, x_3$ are three complex numbers linearly independent over $\Q$ and $y_1, y_2$ are two complex numbers linearly independent over $\Q$ then at least one of the six numbers $\{e^{x_i y_j}: 1 \leq i \leq 3, 1\leq j \leq 2\}$ is transcendental.

Now, suppose by the contrary that $L^\bullet $ is a spiral, that is, $L^\bullet  = \{e^{t(a + \mathfrak{i}b)}: t \in \R\}$ for some $a$ and $b \in \R^*$.
Since $r > 3$ there exist  $\xi_1$, $\xi_2$ and $\xi_3 \in p(\OO^*)$ which are multiplicatively independent over the integers. We may suppose that $\xi_1 = e^{a + \mathfrak{i}}$, $\xi_2 = e^{u(a + \mathfrak{i}b)}$ and $\xi_3 = e^{v(a + \mathfrak{i}b)}$ where $u$ and $v \in \R^*$ and $\mathfrak{i} = \sqrt{-1}$. Remark that $\{1, u, v\}$ are linearly independent over $\Q$, $\{a + \mathfrak{i}b, \mathfrak{i}b\}$ are linearly independent over $\Q$, and the six numbers $\xi_1$, $\xi_2$, $\xi_3$, $\frac{\xi_1}{|\xi_1|}$, $\frac{\xi_2}{|\xi_2|}$, $\frac{\xi_3}{|\xi_3|}$ are all algebraic. This contradicts the six exponentials theorem. \qed

\medskip

If $K_1 = \C$ and $K$ is not a $\mathrm{CM}$-field, it is not difficult to give examples when $L^\bullet $ is the circle group and when $L = \C^*$.

\medskip

\textbf{Examples.}1) For every $n \geq 1$, let $f_n(x) = (x^2 - (\sqrt{n^2 + 1} + n)x + 1)(x^2 - (-\sqrt{n^2 + 1} + n)x + 1)$.
Then $f_n(x)$ is an irreducible polynomial in $\Q[X]$ with two real and two (conjugated) complex roots. Let $K_1 = \Q(\alpha_n)$ where $\alpha_n$ is one of the complex roots of $f_n(x)$. Then $L^\bullet $ is the circle group of $\C^*$.

2) It is easy to see that if $K$ is a totally imaginary, Galois, non-$\mathrm{CM}$-number field of degree $\geq 6$ then $L = \C^*$.

\medskip

Finally, the following is quite plausible:

\medskip

\textbf{Conjecture 2.} $L^\bullet$ is never a spiral.

\medskip

In response to a question of the author, Federico Pellarin observed that Conjecture 2 follows from the still open four exponentials conjecture. Recall that the four exponentials conjecture says that the conclusion of the six exponentials theorem remains valid if replacing the three complex numbers $x_1, x_2, x_3$ by two complex numbers $x_1, x_2$.
The use of the six exponentials theorem in our proof of Proposition \ref{integers1} is inspired by Pellarin's argument
{\footnote{The author is grateful to Federico Pellarin for the useful discussion.}}.

\section{Accumulations points for locally divergent orbits}

As in the introduction, $\Gamma$ is an $\mathcal{S}$-arithmetic subgroup of $G = \G(K_{\mathcal{S}})$ and $T = \te(K_{\mathcal{S}})$ acts on $G/\Gamma$ by left translations.

 In the next lemma $\te$ is identified with $\mathbf{GL}_1^{\rank_K \G}$ via a $K$-isomorphism (see \S 2.2). Under this identification
$\te(\OO)$ is commensurable with $(\OO^*)^{\rank_K \G}$.

\begin{lem}
\label{lem2} Let $h \in \G(K)$. The following assertions hold:
\begin{enumerate}
\item[(a)] There exists a positive integer $m$
such that $\xi\pi(h) = \pi(h)$ for all $\xi \in (\OO^*_m)^{\rank_K \G}$;
\item[(b)] If $h_i$ is a sequence in $G$ such that $\{\pi(h_i)\}$
converges to an element from $G/\Gamma$ then the sequence $\{\pi(h_ih)\}$ admits a converging to an element from $G/\Gamma$ subsequence.
\end{enumerate}
\end{lem}

The lemma is an easy consequence from the commensurability
of $\Gamma$ and $h\Gamma h^{-1}$.

\subsection{Main proposition}

\medskip

Further on we use the notation about the linear algebraic groups as given in \S 2.3.

\begin{prop}\label{main1} Let $n \in \mathcal{N}_{\G}(\te)$ and $\Psi \subset \Pi$. The following conditions are equivalent: 
\begin{enumerate}
\item[(i)] $n \in w_0\mathcal{W}_\Psi$;
\item[(ii)] $\ve_\emptyset^-w_0n\pe_\Psi$ is Zariski dense in $\G$;
\item[(iii)] $w_0n \ve_{\Psi}(w_0n)^{-1} \subset \ve_\emptyset$.
\end{enumerate}
\end{prop}

{\bf Proof.} The implications $\mathrm{(i)} \Rightarrow \mathrm{(ii)}$  and $\mathrm{(i)} \Rightarrow \mathrm{(iii)}$ follow trivially from the definitions in \S 2.3.

Let $\mathrm{(ii)}$ holds. Then $n^{-1}\ve_\emptyset n \pe_\Psi$ is Zariski dense in $\G$. Since $n^{-1}\ve_\emptyset n$ and $\pe_\Psi$ are $\te$-invariant
\begin{equation}
\label{roots1}
\Lie(n^{-1}\ve_\emptyset n) = \Lie(\ve_{\Psi}^-) + \Lie(n^{-1} \ve_\emptyset n \cap \pe_{\Psi}).
\end{equation}
Therefore
\begin{equation}
\label{roots2}
n^{-1}\ve_\emptyset n = \ve_{\Psi}^-  (n^{-1}\ve_\emptyset n \cap \pe_{\Psi}).
\end{equation}
Since $n^{-1}\ve_\emptyset n$ is a product of root groups,
if  $n^{-1}\ve_\emptyset n \cap \ve_\Psi \neq \{e\}$ then, in view of (\ref{roots2}), $n^{-1}\ve_\emptyset n$ contains two opposite root groups which is not possible. So,
$$ n^{-1}\ve_\emptyset n \cap \ve_{\Psi} = \{e\}.$$
This implies $n^{-1}\ve_\emptyset n \cap \pe_{\Psi} = n^{-1}\ve_\emptyset n \cap \mathcal{Z}_\G(\te_{\Psi})$. In view of (\ref{roots1}), $n^{-1}\ve_\emptyset n \cap \mathcal{Z}_\G(\te_{\Psi})$ is a maximal unipotent  subgroup of $ \mathcal{Z}_\G(\te_{\Psi})$. Let  $n' \in \mathcal{W}_\Psi$ be such that $n'(n^{-1}\ve_\emptyset n \cap \mathcal{Z}_\G(\te_{\Psi}))n'^{-1} \subset \ve_\emptyset^-$. Since $n'$ normalizes $\ve_\Psi^-$, it follows from (\ref{roots2}) that
$$
n'n^{-1}\ve_\emptyset nn'^{-1} = \ve_\emptyset^-
$$
which implies $\mathrm{(i)}$.

Suppose that $\mathrm{(iii)}$ holds. Then $ (w_0n)^{-1}\pe_\Psi w_0n \supset \ve_{\Psi}$. Hence, $w_0n \in \pe_\Psi$ by \cite[14.22(iii)]{Borel}. Therefore
$w_0n \in \mathcal{W}_\Psi$, proving $\mathrm{(i)}$. \qed

\medskip

Further on,
 $g = (g_1, g_2, \cdots, g_r) \in G$ where $g_i \in G_i$. We will use the following notational convention: if $h = (h_1, \cdots, h_r) \in G$ and $g_i \in \G(K)$, writing $\pi(hg_i)$ we mean that
 $g_i$ is identified with its the diagonal imbedding in $G$, that is, $hg_i = (h_1g_i, \cdots, h_rg_i) \in G$.

Our main proposition is the following.

\begin{prop}\label{prop1}
Let $\#\mathcal{S} \geq 2$, $\mathrm{rank}_K \mathbf{G} = \mathrm{rank}_{K_{1}} \mathbf{G} = \mathrm{rank}_{K_{2}} \mathbf{G}$, $g_1$ and $g_2 \in \G(K)$ and $\Psi$ be a proper subset of $\Pi$. Let $(s_n, t_n, e, \cdots, e) \in T$ be a sequence and $C > 1$ be a constant such that for all $n$ we have: $|\alpha (s_n)|_{1} > \frac{1}{C}$
 for all $\alpha \in \Pi$, $|\alpha (t_n)|_{2} \underset{n}{\rightarrow} 0$ for all $\alpha \in \Pi \setminus \Psi$ and
 $\frac{1}{C} < |\alpha (t_n)|_{2} < C$ for all $\alpha \in \Psi$. Then the sequence $(s_n, t_n, e, \cdots, e)\pi(g)$ is bounded in $G/\Gamma$ if and only if the
 following conditions are satisfied:
\begin{enumerate}
\item[(i)] $g_1g_2^{-1} \in \ve_{\Psi}^-\pe_{\Psi}$, and
\item[(ii)] there exists a constant $C' > 1$ such that $\frac{1}{C'} < |\alpha (s_n)|_{1}
 \cdot|\alpha (t_n)|_{2} < C'$ for all $\alpha \in \Delta$ and $n$.
\end{enumerate}
\end{prop}

{\bf Proof.} $\Leftarrow)$ Suppose that (i) and (ii) hold. Then $g_1 = vpg_2$ where $v \in \ve_{\Psi}^-(K)$ and $p \in \pe_{\Psi}(K)$. It
follows from the assumption (ii), Lemma \ref{lem2}(a) and Proposition \ref{integers1'}
 that there exists a sequence $d_n \in \textrm{Stab}_T\{\pi(pg_2)\}$ such that the sequence $(s_n d_n^{-1}, t_n d_n^{-1},d_n^{-1}, \cdots, d_n^{-1})$ is bounded in $T$.

Now
\begin{eqnarray*}
& & (s_n, t_n, e, \cdots, e)\cdot(g_1,g_2, \cdots, g_r)\pi(e) =  \\ & & (s_n, t_n, e, \cdots, e)\cdot (vpg_2,p^{-1}pg_2,g_3 \cdots, g_r)\pi(e)= \\ & & (s_nvs_n^{-1}, t_np^{-1}t_n^{-1}, e, \cdots, e)\cdot(s_n,t_n,e, \cdots, e) \cdot \\ && (e,e,g_3(pg_2)^{-1}, \cdots, g_n(pg_2)^{-1})\pi(pg_2)= \\ & &
(s_nvs_n^{-1}, t_np^{-1}t_n^{-1}, e, \cdots, e)\cdot(s_nd_n^{-1},t_nd_n^{-1},e, \cdots, e) \cdot \\ && (e,e,g_3(pg_2)^{-1}d_n^{-1}, \cdots, g_n(pg_2)^{-1}d_n^{-1})\pi(pg_2).
\end{eqnarray*}
Note that $t_n p^{-1}t_n^{-1}$ is bounded in $G_2$. Since $|\alpha (t_n)|_{2} \underset{n}{\rightarrow} 0$ for all $\alpha \in \Pi \setminus \Psi$, it follows from $(ii)$ that $|\alpha (s_n)|_{1} \underset{n}{\rightarrow} \infty$ for all $\alpha \in \Pi \setminus \Psi$. Therefore $s_n v s_n^{-1} \underset{n}{\rightarrow} e$ in $G_1$.
 Now using the choice of $d_n$ we conclude that
$(s_n, t_n, e, \cdots, e)\pi(g)$ is bounded in $G/\Gamma$.

$\Rightarrow)$ Let $(s_n, t_n, e, \cdots, e)\pi(g)$ be bounded. Using Bruhat decomposition one can write $g_1 g_2^{-1} = v^-w_0np$ where $v^- \in \ve_\emptyset^-$, $n \in \mathcal{N}_\G(\te)$,
$p \in \pe_\Psi$ and, as usual, $w_0\ve_\emptyset^-w_0^{-1} = \ve_\emptyset$. Suppose that $g_1 g_2^{-1} \notin \ve_\Psi^-\pe_\Psi$. Since $\ve_\emptyset^-\pe_\Psi = \ve_\Psi^-\pe_\Psi$ it follows from Proposition \ref{main1} the
existence of a root $\chi$ such that $\mathbf{U}_\chi \subset \ve_\Psi$ and $\chi \circ \textrm{Int}(w_0n)^{-1}$ is a negative root.
As above, using Lemma \ref{lem2}(a) and Proposition \ref{integers1'}, we fix a sequence
$d_n \in T \cap \textrm{Stab}_G\{\pi(pg_2)\}$ such that the sequence $\{t_nd_n^{-1}\}$ is bounded in $G_2$ and the sequence $\{d_n\}$ is bounded in every $G_i$, $i \geq 3$. Since $|\alpha (s_n)|_{1} > \frac{1}{C}$
 for all $\alpha \in \Pi$ we get that $\{|\chi((w_0n)^{-1}s_nw_0n)|_{1}\}$ is bounded from above and since $|\alpha (t_n)|_{2} \underset{n}{\rightarrow} 0$ for all $\alpha \in \Pi \setminus \Psi$ it follows from Proposition \ref{integers1'} and the choice of $\{d_n\}$ that
$|\chi(d_n^{-1})|_{1} \underset{n}{\rightarrow} 0$.
Therefore
\begin{equation}
\label{pf.main}
|\chi((w_0n)^{-1}s_nw_0n)\chi(d_n^{-1})|_{1} \underset{n}{\rightarrow} 0.
\end{equation}
Hence $\{(w_0n)^{-1}s_n(w_0n)d_n^{-1}\}$ is unbounded in $T_1$.

Note that
\begin{eqnarray*}
& & (s_n, t_n, e, \cdots, e)\pi(g) =  \\ & & (s_n, t_n, e, \cdots, e)\cdot (v^{-}w_0n,p^{-1},g_3(pg_2)^{-1},\cdots,g_r(pg_2)^{-1})\pi(pg_2)= \\ & &
((s_nv^-s_n^{-1}w_0n)((w_0n)^{-1}s_n(w_0n)d_n^{-1}), (t_np^{-1}t_n^{-1})t_nd_n^{-1}, g_3(pg_2)^{-1}d_n^{-1},\\ & & \cdots, g_r(pg_2)^{-1}d_n^{-1})\pi(pg_2).
\end{eqnarray*}
where $\{s_nv^-s_n^{-1}w_0n\}$ is bounded in $G_1$, $\{t_np^{-1}t_n^{-1})t_nd_n^{-1}\}$ is bounded in $G_2$ and $\{g_i(pg_2)^{-1}d_n^{-1}\}$ is bounded in $G_i$ for every $i \geq 3$. Since $T_1\pi(pg_2)$ is divergent, it follows from the above that $(s_n, t_n, e, \cdots, e)\pi(g)$ is unbounded, contradicting our hypothesis. We have proved (i).

Let $g_1 = v^-pg_2$, where $v^- \in \ve_{\Psi}^-, \ p\in \pe_\Psi$. Then, with $d_n$ chosen as above, we have
\begin{eqnarray*}
& & (s_n, t_n, e, \cdots, e)\pi(g) = ((s_nv^-s_n^{-1})(s_nd_n^{-1}),
(t_np^{-1}t_n^{-1})(t_nd_n^{-1}), \\ & & g_3(pg_2)^{-1}d_n^{-1},\cdots, g_r(pg_2)^{-1}d_n^{-1})\pi(pg_2).
\end{eqnarray*}
Since $\{s_nv^-s_n^{-1}\}$ is bounded in $G_1$, $\{t_np^{-1}t_n^{-1}\}$ and $\{t_nd_n^{-1}\}$ are both bounded in $G_2$, and
$\{d_n\}$ is bounded in $G_i$ for each $i \geq 3$ it follows from the assumptions that $(s_n, t_n, e, \cdots, e)\pi(g)$ is bounded in $G/\Gamma$
and $T_1\pi(g)$ is divergent that $\{s_nd_n^{-1}\}$ is bounded in $G_1$. Hence there exists $C_1 > 1$ such that $\frac{1}{C_1} < |\alpha(s_n d_n^{-1})|_{1}\cdot |\alpha(t_n d_n^{-1})|_{2} < C_1$ for all $\alpha \in \Pi$. By Artin's product formula $\prod_{v \in \mathcal{V}}|\alpha(d_n)|_{v} = 1$ where $\mathcal{V}$ is the set of all normalized valuations of $K$. This implies (ii).  \qed

The above proposition implies:

\begin{cor}
\label{THM1}
Let $s_n \in T_1$ and $t_n \in T_2$ be such that for every $\alpha \in \Phi$ each of the sequences $|\alpha(s_n)|_1$ and $|\alpha(t_n)|_2$ converges to an element from $\mathbb{R} \cup \infty$. We suppose that $g_1$ and $g_2 \in \G(K)$ and that $(s_n, t_n, e, \cdots, e)\pi(g)$
converges in $G/\Gamma$. Then there exist $\Psi \subset \Pi$ and $\omega_1, \omega_2 \in \mathcal{W}$ with the following properties:
\begin{enumerate}
\item[(i)] $\omega_1\mathbf{P}_\Psi^-\omega_1^{-1}(K_1) \times \omega_2\mathbf{P}_\Psi\omega_2^{-1}(K_2)= \\
  = \{(x,y) \in {G}_1 \times {G}_2 : \mathrm{Int}(s_n,t_n)(x,y) \ {is} \ {bounded} \ {in} \ {G}_1 \times {G}_2\},$
\item[(ii)] $g_1g_2^{-1} \in \omega_1 \ve_{\Psi}^-\pe_{\Psi}\omega_2^{-1}$.
\item[(iii)] If $g_1g_2^{-1} = \omega_1 v_{\Psi}^-z_\Psi v_{\Psi}\omega_2^{-1}$,
where $v_{\Psi}^- \in \ve_{\Psi}^-(K)$,
$z_\Psi \in \mathcal{Z}_\G(\te_\Psi)(K)$ and $v_{\Psi} \in \ve_{\Psi}(K)$, then
$$\underset{n}{\lim} ({s}_n, {t}_n, e, \cdots, e)\pi({g}) \in T(\omega_1(v_{\Psi}^-)^{-1}\omega_1^{-1}, \omega_2 v_{\Psi}\omega_2^{-1}, e, \cdots, e)\pi(g).$$
\end{enumerate}
\end{cor}
{\bf Proof.} Since $|\alpha(s_n)|_1$ converges for every $\alpha \in \Phi$, there exists a parabolic $K$-subgroup $\mathbf{P}$ containing $\mathbf{T}$ such that
$\mathbf{P}(K_1) = \{x \in \mathbf{G}(K_1): \mathrm{Int}(s_n)x \ \mathrm{is} \ \mathrm{bounded}\}$. Let $\Psi \subset \Pi$ and $\omega_1 \in \mathcal{W}$ be such that $\mathbf{P} = \omega_1 \mathbf{P}_\Psi^- \omega_1^{-1}$.
Similarly, we find $\Psi' \subset \Pi$ and $\omega_2 \in \mathcal{W}$ such that
$\omega_2\mathbf{P}_{\Psi'}\omega_2^{-1}(K_2) = \{x \in \mathbf{G}(K_2): \mathrm{Int}(t_n)x \ \mathrm{is} \ \mathrm{bounded}\}$.

Put $\widetilde{g} = (\omega_1^{-1}g_1, \omega_2^{-1}g_2, g_3, \cdots, g_r)$, $\widetilde{s}_n = \omega_1^{-1}s_n\omega_1$ and $\widetilde{t}_n  = \omega_2^{-1}t_n\omega_2$. Then $(\widetilde{s}_n, \widetilde{t}_n, e, \cdots, e)\pi(\widetilde{g})$ converges, $\mathbf{P}_\Psi^{-}(K_1) = \{x \in \mathbf{G}(K_1): \mathrm{Int}(\widetilde{s}_n)x$ $\mathrm{is} \ \mathrm{bounded}\}$ and $\mathbf{P}_{\Psi'}(K_2) = \{x \in \mathbf{G}(K_2): \mathrm{Int}(\widetilde{t}_n)x \ \mathrm{is} \ \mathrm{bounded}\}$. By the choice of $\Psi$ and $\Psi'$, there exists a constant $C > 1$ such that $C > |\alpha(\widetilde{s}_n)|_1 > \frac{1}{C}$ for all $\alpha \in \Psi$, $|\alpha(\widetilde{s}_n)|_1 \underset{n}{\rightarrow} \infty$ for all $\Pi \setminus \Psi$, $|\alpha(\widetilde{t}_n)|_2 \underset{n}{\rightarrow} 0$ for all $\Pi \setminus \Psi'$ and $C > |\alpha(\widetilde{t}_n)|_2 > \frac{1}{C}$ for all $\alpha \in \Psi'$. Replacing, if necessary,
$C$ by a larger constant we may suppose that for every $\alpha \in \Phi$ either $C > |\alpha(\widetilde{t}_n)|_2 > \frac{1}{C}$ for all $n$
  or $|\alpha(\widetilde{t}_n)|_2$ is converging to $0$ or $\infty$. It follows from Proposition \ref{prop1}(ii) that
$\frac{1}{C} < |\alpha (\widetilde{s}_n)|_{1} \cdot|\alpha (\widetilde{t}_n)|_{2} < C$ for all $\alpha \in \Pi$ and $n$. This implies easily
that $\Psi = \Psi'$. In view of Proposition \ref{prop1}(i) $g_1g_2^{-1} \in \omega_1 \ve_{\Psi}^-\pe_{\Psi}\omega_2^{-1}$. Hence (i) and (ii) hold.

Let $\omega_1^{-1}g_1 =  v_{\Psi}^-z_\Psi v_{\Psi}\omega_2^{-1} g_2$. Then
\begin{eqnarray*}
& &
(\widetilde{s}_n, \widetilde{t}_n, e, \cdots, e)\pi(\widetilde{g}) = \\ & &
(\widetilde{s}_n, \widetilde{t}_n, e, \cdots, e)(v_{\Psi}^-z_\Psi, v_{\Psi}^{-1}, g_3(v_{\Psi}\omega_2^{-1}g_2)^{-1}, \cdots, g_r(v_{\Psi}\omega_2^{-1}g_2)^{-1})\pi(v_{\Psi}\omega_2^{-1}g_2),
\end{eqnarray*}
Lemma \ref{lem2}(a) and Proposition \ref{integers1'} imply the existence of $d_n \in \te_\Psi \cap \textrm{Stab}_G\{\pi(v_{\Psi}\omega_2^{-1}g_2)\}$ such that $\{\widetilde{s}_n d_n^{-1}\}$ is bounded in $G_1$, $\{\widetilde{t}_n d_n^{-1}\}$ is bounded in $G_2$ and $\{d_n\}$ is bounded in $G_i$ for all $i \geq 3$. Since $d_n v_{\Psi}^- d_n^{-1} \underset{n}{\rightarrow} 0$ in $G_1$ and $d_n v_{\Psi}^{-1} d_n^{-1} \underset{n}{\rightarrow} 0$ in $G_2$, we get that
$$
\underset{n}{\lim} (\widetilde{s}_n, \widetilde{t}_n, e, \cdots, e)\pi(\widetilde{g}) \in
T((v_{\Psi}^-)^{-1}\omega_1^{-1},  v_{\Psi}\omega_2^{-1}, e, \cdots, e)\pi(g),
$$
which implies (iii).
\qed

\medskip

\section{Closures of locally divergent orbits for $\#\mathcal{S} = 2$}

In this section $\mathcal{S} = \{v_1, v_2\}$, $g = (g_{1},g_{1}) \in G$
and $T\pi(g)$ is a locally divergent orbit.

We continue to use the notation $\Pi$, $\Psi \subset \Pi$, $\mathbf{P}_\Psi$, $\mathbf{P}^-_\Psi$, $\ve_\Psi$,
$\mathbf{V}_\Psi^-$ and $\mathcal{Z}_{\G}(\te_\Psi)$ as introduced in \S 2.3. Further on by a parabolic subgroup we
mean a parabolic subgroup defined over $K$.

For every $\Psi \subset \Pi$, we put
$$\mathcal{P}_\Psi(g) = \{\omega_1{\mathbf{P}}_\Psi^- \omega_1^{-1} \times \omega_2\mathbf{P}_\Psi \omega_2^{-1}| \omega_1, \omega_2
\in \mathcal{N}_\G(\te), g_{1}g_{2}^{-1} \in \omega_1\mathbf{V}_\Psi^-\mathbf{P}_\Psi \omega_2^{-1}\} \footnote{Note that $\pe^-_\Psi$ and $\pe_\Psi$ are not always conjugated. Therefore we can not replace $\pe^-_\Psi$ by $\pe_\Psi$ in the definition of $\mathcal{P}_\Psi(g)$ . For exemple, if $\G$ is of type $D_l$, $l \geq 4,$ and $\alpha \in \Pi$ be such that $\omega_0(\alpha) \neq -\alpha$ then $\pe^-_{\{\alpha\}}$ is not conjugated to $\pe_{\{\alpha\}}$.}.$$
It is trivial but worth mentioning that $\mathcal{P}_\Psi(g)$ is a \textit{finite} set of parabolic subgroups of $\G \times \G$ and $\omega_1{\mathbf{P}}_\Psi^- \omega_1^{-1} \times \omega_2\mathbf{P}_\Psi \omega_2^{-1} \in \mathcal{P}_\Psi(g)$ if and only if $g_{1}g_{2}^{-1}$ belongs to the \textit{non-empty Zariski open subset} $\omega_1\mathbf{V}_\Psi^-\mathbf{P}_\Psi \omega_2^{-1}$. It is clear that $\mathcal{P}_\Pi(g) = \{\G \times \G\}$ and $\mathcal{P}_\emptyset(g)$ consists of minimal parabolic
$K$-subgroups of $\G \times \G$.

Denote
$$\mathcal{P}(g) = \underset{\Psi \subset \Pi}{\bigcup}\mathcal{P}_\Psi(g).$$

To every $\mathbf{P} \in \mathcal{P}(g)$ we associate  a locally divergent $T$-orbit as follows.
Let $\mathbf{P} = \omega_1{\mathbf{P}}_\Psi^- \omega_1^{-1} \times \omega_2\mathbf{P}_\Psi \omega_2^{-1}$ and
\begin{equation}
\label{decomposition}
g_{1}g_{2}^{-1} = \omega_1 v_{\Psi}^-z_\Psi v_{\Psi}\omega_2^{-1},
\end{equation}
where $v_{\Psi}^- \in \mathbf{V}_\Psi^-$,
$z_\Psi \in \mathcal{Z}_{\G}(\te_\Psi)$, $v_{\Psi} \in \mathbf{V}_\Psi$ $\omega_1$ and $\omega_2 \in \mathcal{N}_\G(\te)$. The locally divergent orbit associated to $\mathbf{P}$ is
\begin{equation}
\label{explicit}
\mathrm{Orb}_{g}(\pe) \bydefn T(\omega_1(v_{\Psi}^-)^{-1}\omega_1^{-1}, \omega_2v_{\Psi}\omega_2^{-1})\pi(g).
\end{equation}
Clearly, $\mathrm{Orb}_{g}(\G \times \G) = T\pi(g)$.

Concerning (\ref{explicit}), note that the matrices $\omega_1, \omega_2, v_{\Psi}^-$ and $v_{\Psi}$ are with coefficients
from the universal domain $\hat{K}$ (see 2.2) and they are not uniquely
defined by the decomposition (\ref{decomposition}). Proposition \ref{inv}(a) shows that the orbit
$\mathrm{Orb}_{g}(\pe)$ is well-defined, i.e., it does not depend on the decomposition (\ref{decomposition}).

\begin{prop}
\label{inv}
Let $\mathbf{P} \in \mathcal{P}(g)$.
Then
\begin{enumerate}
\item[(a)] $\mathrm{Orb}_{g}(\pe)$ is well-defined and locally divergent;
\item[(b)] If $w \in \mathcal{W} \times \mathcal{W}$ then $w\pe w^{-1} \in \mathcal{P}(wg)$ and
\begin{equation}
\label{invariance}
w\mathrm{Orb}_{g}(\pe) = \mathrm{Orb}_{wg}(w\pe w^{-1}).
\end{equation}
\end{enumerate}
\end{prop}
{\bf Proof.} (a) The decomposition (\ref{decomposition}) is determined by the choice of $\omega_1$ and  $\omega_2 \in \mathcal{N}_\G(\te)$. If
$\omega_1 = a_1 \widetilde{\omega}_1$ and $\omega_2 = a_2 \widetilde{\omega}_2$, where $a_1$ and  $a_2 \in \mathcal{Z}_\G(\te)$, then $g_{1}g_{2}^{-1} = \omega_1 v_{\Psi}^-z_\Psi v_{\Psi}\omega_2^{-1} = \widetilde{\omega}_1 \widetilde{v}_{\Psi}^-\widetilde{z}_\Psi \widetilde{v}_{\Psi}\widetilde{{\omega}}_2^{-1}$, where $\widetilde{v}_{\Psi}^- = \mathrm{Int}(\widetilde{{\omega}}_1^{-1}a_1\widetilde{{\omega}}_1)({v}_{\Psi}^-)$, $\widetilde{z}_\Psi = (\widetilde{{\omega}}_1^{-1}a_1\widetilde{{\omega}}_1){z}_\Psi(\widetilde{{\omega}}_2^{-1}a_2\widetilde{{\omega}}_2)^{-1}$ and $\widetilde{v}_{\Psi} = \mathrm{Int}(\widetilde{{\omega}}_2^{-1}a_2\widetilde{{\omega}}_2)^{-1}({v}_{\Psi})$. Therefore $\omega_1(v_{\Psi}^-)^{-1}\omega_1^{-1} = \widetilde{\omega}_1(\widetilde{v}_{\Psi}^-)^{-1}\widetilde{\omega}_1^{-1}$
and $\omega_2 v_{\Psi}\omega_2^{-1} = \widetilde{\omega}_2\widetilde{v}_{\Psi}\widetilde{\omega}_2^{-1}$, proving that (\ref{explicit}) does not depend on
the choice of $\omega_1$ and  $\omega_2$.

In remains to prove that $\mathrm{Orb}_{g}(\pe)$ is locally divergent. Since $T\pi(g)$ is locally divergent, there exist
$a = (a_1, a_2) \in \mathcal{Z}_{G}(T)$ and $\widetilde{g} = (\widetilde{{g}}_1, \widetilde{{g}}_2) \in \G(K) \times \G(K)$ such that $g = a \widetilde{g}$.
It is well-known that $\mathcal{N}_\G(\te) = \mathcal{Z}_\G(\te) \mathcal{N}_\G(\te)(K)$ (\cite[Theorem 21.2]{Borel}). We choose $\omega_1 = a_1 \widetilde{\omega}_1$ and $\omega_2 = a_2 \widetilde{\omega}_2$ with $\widetilde{\omega}_1$ and  $\widetilde{\omega}_2 \in \mathcal{N}_\G(\te)(K)$. With the above notation, $g_{1}g_{2}^{-1} = a_1\widetilde{g}_{1}\widetilde{g}_{2}^{-1}a_2^{-1} = a_1 \widetilde{\omega}_1 v_{\Psi}^-z_\Psi v_{\Psi} \widetilde{\omega}_2^{-1}{a}_2^{-1}$. So, $\widetilde{\omega}_1^{-1}\widetilde{g}_{1}\widetilde{g}_{2}^{-1}\widetilde{\omega}_2 = v_{\Psi}^-z_\Psi v_{\Psi} \in \G(K)$. Since the product map $\mathbf{V}_\Psi^- \times \mathcal{Z}_{\G}(\te_\Psi) \times \mathbf{V}_\Psi \rightarrow $ $\mathbf{V}_\Psi^- \mathcal{Z}_{\G}(\te_\Psi)\mathbf{V}_\Psi$  is a regular $K$-isomorphism (cf. \cite[Theorem 21.15]{Borel}), we get that $v_{\Psi}^-, z_\Psi$ and $v_{\Psi} \in \G(K)$. Therefore
$$
\mathrm{Orb}_{g}(\pe) = aT(\widetilde{\omega}_1(v_{\Psi}^-)^{-1}\widetilde{\omega}_1^{-1}, \widetilde{\omega}_2v_{\Psi}\widetilde{\omega}_2^{-1})\pi(\widetilde{g})
$$
completing the proof.

The part (b) follows from the definition (\ref{explicit}) by a simple computation.
\qed

\begin{thm}
\label{THM2} With the above notation, let $\pe \in \mathcal{P}(g)$. Then
\begin{equation*}
\label{closure}
\overline{\mathrm{Orb}_{g}(\pe)} = \bigcup_{\pe' \in \mathcal{P}(g), \ \pe' \subset \pe} \mathrm{Orb}_{g}(\pe').
\end{equation*}
In particular,
$$
\overline{T\pi(g)} = \bigcup_{\pe \in \mathcal{P}(g)} \mathrm{Orb}_{g}(\pe).
$$
\end{thm}

Theorem \ref{THM2} implies immediately:

\begin{cor}
\label{cor1} Let ${\pe} \in \mathcal{P}(g)$. Choose a $\widetilde{\pe} \in \mathcal{P}(g)$ with ${\mathrm{Orb}_{g}(\widetilde{\pe})} = {\mathrm{Orb}_{g}({\pe})}$ and being minimal with this property.
Then
\begin{eqnarray*}
\label{closure1}
&& \overline{\mathrm{Orb}_{g}(\pe)} \setminus \mathrm{Orb}_{g}(\pe) =
\bigcup_{\pe' \in \mathcal{P}(g), \ \pe' \subsetneqq {\widetilde{\pe}}} \overline{\mathrm{Orb}_{g}(\pe')}.
\end{eqnarray*}
In particular, $T\pi(g)$ is open in its closure.
\end{cor}

\subsection{Proof of Theorem \ref{THM2}}
In view of  (\ref{invariance}), it is enough to prove the theorem for $\pe = \pe_\Psi^- \times \pe_\Psi$.
In this case $g_1 g_2^{-1} = v_\Psi^- z_\Psi v_\Psi$, where $v_\Psi^- \in \ve_\Psi^-$, $z_\Psi \in \mathcal{Z}_{\G}(\te_\Psi)$
and $v_\Psi \in \ve_\Psi$. It follows from (\ref{explicit}) that
 \begin{equation}
\label{formula0}
 \mathrm{Orb}_g(\pe) = T(z_\Psi,e)\pi(v_\Psi g_2).
 \end{equation}
Note that $\mathcal{Z}_{\G}(\te_\Psi)$ is a reductive $K$-algebraic group and the orbit $\mathcal{Z}_{G}(T_\Psi)\pi(v_\Psi g_2)$ is closed containing  $\overline{\mathrm{Orb}_g(\pe)}$.
The $T$-orbits on $\mathcal{Z}_{G}(T_\Psi)\pi(v_\Psi g_2)$ contained in $\overline{\mathrm{Orb}_g(\pe)}$ are given by Corollary \ref{THM1} (iii) applied to $\mathcal{Z}_{\G}(\te_\Psi)$ instead of $\G$. More precisely, if $Tm$ is such an orbit there exist $\Psi' \subset \Psi$ and $(\omega_1, \omega_2) \in \mathcal{W}_{\Psi} \times \mathcal{W}_{\Psi}$ such that $\pe' \bydefn \omega_1{\mathbf{P}}_{\Psi'}^- \omega_1^{-1} \times \omega_2\mathbf{P}_{\Psi'} \omega_2^{-1} \in \mathcal{P}(g)$ and $Tm = \mathrm{Orb}_g(\pe')$. Hence $\overline{\mathrm{Orb}_{g}(\pe)} \subset \bigcup_{\pe' \in \mathcal{P}(g), \ \pe' \subset \pe} \mathrm{Orb}_{g}(\pe')$. In order to prove the opposite inclusion, choose a sequence $(s_n, t_n) \in T_1 \times T_2$ such that $\alpha(\omega_1^{-1}s_n \omega_1) = \alpha(\omega_2^{-1}t_n \omega_2) = 1$ for all $\alpha \in \Psi'$, $|\alpha(\omega_1^{-1}s_n \omega_1)|_1 \underset{n}{\rightarrow} \infty$ for all $\alpha \in \Pi \setminus \Psi'$ and $|\alpha(\omega_2^{-1}t_n \omega_2)|_2 \underset{n}{\rightarrow} 0$ for all $\alpha \in \Pi \setminus \Psi'$. It is easy to see that $(s_n, t_n)(z_\Psi,e)\pi(v_\Psi g_2)$ converges to an element from $\mathrm{Orb}_{g}(\pe')$
 completing the proof of the theorem. \qed

\medskip

The theorem implies that the closed $T$-orbits in $\overline{T\pi(g)}$ are parameterized by the elements of $\mathcal{P}_\emptyset(g)$, that is, by the minimal parabolic subgroups of $\G \times \G$ belonging to $\mathcal{P}(g)$.

\begin{cor}
\label{cor2} If $\pe$ is minimal in $\mathcal{P}(g)$ then $\mathrm{Orb}_{g}(\pe)$ is closed and $\pe$ is a minimal parabolic subgroup of $\G \times \G$. In particular, $\mathcal{P}_\emptyset(g) \neq \emptyset$ and $\{\mathrm{Orb}_{g}(\pe): \pe \in \mathcal{P}_\emptyset(g)\}$ is the set of all closed $T$-orbits in $\overline{T\pi(g)}$.
\end{cor}
{\bf Proof.}
 By (\ref{invariance}) and Theorem \ref{ldo}(a) we may (and will) suppose that $g_1$ and $g_2 \in \G(K)$ and
 $\pe = {\mathbf{P}}_\Psi^-  \times \mathbf{P}_\Psi$. In this case ${\mathrm{Orb}}_{g}(\pe)$ is explicitly given by the formula (\ref{formula0}).
If $\pe$ is minimal among the subgroups in $\mathcal{P}(g)$ then $\mathrm{Orb}_{g}(\pe)$ is closed  in view of Theorem \ref{THM2}.
It follows from Theorem \ref{ldo} (b) that $z_\Psi \in \mathcal{N}_{Z_{\G}(\te_\Psi})(\te)$. So,
$$
g_1g_2^{-1} = v_{\Psi}^-(z_{\Psi} v_{\Psi}z_{\Psi}^{-1})z_{\Psi} \in \ve_{\emptyset}^{-}\pe_{\emptyset}z_{\Psi},
$$
  implying that $\pe_{\emptyset}^{-} \times z_{\Psi}^{-1}\pe_{\emptyset} z_{\Psi} \in \mathcal{P}(g)$. Since $\pe$ is minimal in  $\mathcal{P}(g)$ and $\pe_{\emptyset}^{-} \times z_{\Psi}^{-1}\pe_{\emptyset} z_{\Psi} \subset {\mathbf{P}}_\Psi^-  \times \mathbf{P}_\Psi$, we get that $\pe = \pe_{\emptyset}^{-} \times z_{\Psi}^{-1}\pe_{\emptyset} z_{\Psi}$, i.e. $\pe$ is a minimal parabolic subgroup of $\G \times \G$. \qed

Theorem \ref{THM2} easily implies the following refinement of Theorem \ref{ldo}(b) for $\#\mathcal{S} = 2$.

\begin{cor}
\label{cor4}
The following conditions are equivalent:
\begin{enumerate}
\item[(a)] ${T\pi(g)}$ is closed,
\item[(b)] $\overline{T\pi(g)}$ is homogenous,
\item[(c)] $g \in \mathcal{N}_G(T)\G(K)$,
\end{enumerate}
\end{cor}
{\bf Proof.} In view of Theorem \ref{ldo}(b), we need only to prove that $(b) \Rightarrow (a)$.
If $\overline{T\pi(g)}$ is homogeneous then $\overline{T\pi(g)} = H\pi(g)$, where $H$ is a closed subgroup of $G$ containing $T$. Since $\overline{T\pi(g)}$ is a finite union
of $T$-orbits, $T$ is a subgroup of finite index in $H$. Therefore ${T\pi(g)}$ is closed.
\qed

\medskip

It is easy to see that the map $\pe \mapsto {\mathrm{Orb}_{g}(\pe)}$, $\pe \in \mathcal{P}(g)$, is not always injective.
It becomes injective if $g_{1}g_{2}^{-1}$ belongs to a non-empty Zariski dense subset of $\G$.

\begin{cor}
\label{cor3}
For every $\Psi \subset \Pi$, denote by $n_{\Psi}$ the number of parabolic subgroups containing $\T$ and conjugated to $\mathbf{P}_{\Psi}$. We have
\begin{enumerate}
\item[(a)] The number of different $T$-orbits in $\overline{T\pi(g)}$ is bounded from above by $\sum_{\Psi \subset \Pi}n_{\Psi}^2$ and the number of different closed $T$-orbits in $\overline{T\pi(g)}$ is bounded from above by $n_{\emptyset}^2$;
\item[(b)] Given $g_2 \in \G(K)$, there exists a non-empty Zariski dense subset $\Omega \subset \G(K)$ such that if $g_{1}g_{2}^{-1} \in \Omega$
then
$\overline{T\pi(g)}$ is a union of exactly $\sum_{\Psi \subset \Pi}n_{\Psi}^2$ pairwise different $T$-orbits and among them exactly $n_{\emptyset}^2$ are
closed. In particular, the map $\mathrm{Orb}_{g}(\cdot)$ is injective.
\end{enumerate}
\end{cor}

{\bf Proof.} The part (a) is an immediate consequence from Theorem \ref{THM2} and the definition of $\mathrm{Orb}_{g}(\pe)$.

Let us prove (b). Denote by $\mathcal{P}$ the set of all parabolic subgroups $$\omega_1\pe_{\Psi}^-\omega_1^{-1} \times \omega_2\pe_{\Psi}\omega_2^{-1}$$
where $\omega_1$ and $\omega_1 \in \mathcal{N}_\G(\T)$ and $\Psi \subset \Pi$. Let $\Omega_1 = \underset{(\omega_1,\omega_1) \in \mathcal{W} \times \mathcal{W}}{\bigcap} \omega_1^{-1}
\pe_{\emptyset}^-\pe_{\emptyset}\omega_2$. Then $\Omega_1$ is $\mathcal{W}$-invariant, Zariski open,
non-empty and $\mathcal{P} = \mathcal{P}(g)$ if and only if $g_1g_2^{-1} \in \Omega_1$.

Since every parabolic subgroup of $\G \times \G$ containing $\te \times \te$ is generated by its minimal parabolic subgroups containing $\te \times \te$,
it is enough to prove the existence of a Zariski dense $\Omega \subset \G(K) \cap \Omega_1$ such that the restriction of $\mathrm{Orb}_{g}(\cdot)$ to
the set of minimal parabolic subgroups of $\G \times \G$ containing $\te \times \te$ is injective whenever $g_1g_2^{-1} \in \Omega$.

Let $\mathcal{W}_\circ \subset \mathcal{N}_{\G}(\te)(K)$ be a finite set containing $e$ such that the natural projection
$\mathcal{W}_\circ \rightarrow \mathcal{W}$
is bijective.

Let $\Delta = g_2 \Gamma g_2^{-1}$. Since the product map
$\ve_{\emptyset}^- \times \mathcal{Z}_\G(\mathbf{T}) \times \ve_{\emptyset} \rightarrow \G$ is a $K$-rational
isomorphism, the projection $\ve_{\emptyset}^- \times \mathcal{Z}_\G(\mathbf{T})  \times \ve_{\emptyset} \rightarrow \ve_{\emptyset}^-$ induces a rational
map $p: \G \rightarrow \ve_{\emptyset}^-$ whose restriction on $\Omega_1$ is regular.
Fix a non-archimedean completion $F$ of $K$ different
from $K_1$ and $K_2$. Let $\overline{p(\Delta)}$ be the closure of $p(\Delta)$ in
$\mathbf{V}_\emptyset^-(F)$ for the topology on $\mathbf{V}_\emptyset^-(F)$ induced by the topology on $F$. Then $\overline{p(\Delta)}$ is compact
in $\mathbf{V}_\emptyset^-(F)$.
There exists a non-empty Zariski open subset $\Omega_2 \subset \Omega_1$ such that if $x \in \Omega_2 \cap \G(F)$,
$\omega \in \mathcal{W}_\circ \setminus \{e\}$ and $\omega x \omega^{-1} =v^-zv$,
where $v^- \in \ve_{\emptyset}^-(F)$, $z \in \mathcal{Z}_\G(\mathbf{T})(F)$ and $v \in \ve_{\emptyset}(F)$,
then $p(\omega^{-1}v ^{-1}\omega) \neq e$. Moreover, there exists a compact $C \subset \G(F)$ such that if, with the
above notation, $x \in \Omega_2 \cap \G(F) \setminus C$ then $p(\omega^{-1}v ^{-1}\omega) \notin \overline{p(\Delta)}$.

Let $\Omega = \Omega_2 \cap \G(K) \setminus C$. It is clear that $\Omega$ is non-empty and Zariski dense in $\G$.
Let $g_1 \in \Omega g_2$. Let $\pe$ and $\pe' \in \mathcal{P}_\emptyset$ be such that $\mathrm{Orb}_g(\pe) = \mathrm{Orb}_g(\pe')$.
We need to prove that $\pe = \pe'$.
In view of (\ref{invariance}), we may assume that
$\pe = \pe_{\emptyset}^- \times \pe_{\emptyset}$ and  $\pe' = \omega_1^{-1}\pe_{\emptyset}^-\omega_1 \times \omega_2^{-1}\pe_{\emptyset}\omega_2$,
where $\omega_1$ and $\omega_1 \in \mathcal{W}_\circ$. Then
$$
g_1 g_2^{-1} = v_1^-z_1 v_1 = \omega_1^{-1} v^-z v\omega_2,
$$
where $v^-$ and $v_{1}^- \in \ve_{\emptyset}^-(K)$, $z$ and
$z_1 \in \mathcal{Z}_\G(\T)(K)$, and $v$ and $v_{1} \in \ve_{\emptyset}(K)$. Using (\ref{explicit}) we get $(t_1,t_2) \in T$
and $\delta \in \Delta$ such that
\begin{equation}
\label{formula}
t_1z_1v_1 = \omega_1^{-1}z v\omega_2 \delta \ \mathrm{and} \ t_2v_1 = \omega_1^{-1} v\omega_2 \delta.
\end{equation}
This implies
$$
\omega_1^{-1}\omega_2 = (t_1 z_1 t_2^{-1})(\omega_2^{-1}z \omega_2) \in \mathcal{Z}_\G(\T).
$$
Therefore $\omega_1 = \omega_2 = \omega$.
Using (\ref{formula})
$$
(\omega^{-1}v^{-1}\omega)t_2 v_1 \in \Delta.
$$
Hence
$$p(\omega^{-1}v^{-1}\omega) \in \overline{p(\Delta)}$$
implying that $\omega = e$, i.e. $\pe = \pe'$.
\qed

\section{Proofs of Theorem \ref{r>2} and Theorem \ref{cor5}}

We assume that $r = \# \mathcal{S} > 2$ and $\rank_K \G = \rank_{K_v} \G > 0$ for all $v \in \mathcal{S}$. Recall that
$\te$ is a maximal $K$-split torus in $\G$, $G_v = \G(K_v)$
and $T_v = \te(K_v)$. We let $\mathcal{S} = \{v_1, \cdots, v_r\}$ and use often the simpler notation $K_i = K_{v_i}$,
$G_i = G_{v_i}$ and $T_i = T_{v_i}$.

\subsection{Horospherical subgroups.}

Let $t \in T_v, v \in \mathcal{S}$. We set
$$
W^+(t) = \{x \in G_v: \underset{n \rightarrow +\infty}{\lim} \ t^{-n}xt^n = e\},
$$

$$
W^-(t) = \{x \in G_v: \underset{n \rightarrow +\infty}{\lim} \ t^{n}xt^{-n} = e\}
$$
and
$$
Z(t) = \{x \in G_v: t^{n}xt^{-n}, \ n \in \Z, \ \mathrm{is} \ \mathrm{bounded}\}.
$$
Then $W^+(t)$ (respectively, $W^-(t)$) is the positive (respectively, negative) horospherical subgroup of $G_v$ corresponding to
$t$.

The next proposition is well known. It follows easily from the assumption  that $\rank_K \G = \rank_{K_v} \G > 0$ for all $v \in \mathcal{S}$.

\begin{prop}\label{horospherical} With the above notation, there exist a basis $\Pi$ of $\Phi(\te, \G)$ and $\Psi \subset \Pi$
such that $Z(t) = \mathcal{Z}_{\G}({\te_\Psi})(K_v)$, $W^+(t) = \mathcal{R}_u(\mathbf{P}_\Psi)(K_v)$, and
 $W^-(t) = \mathcal{R}_u(\mathbf{P}_\Psi^-)(K_v)$.
\end{prop}

\begin{lem}\label{units} Let $\Psi \subset \Pi$, $\sigma \in \G(K)$ and $1 \leq s_1 < s_2 \leq r$.
There exists a sequence $t_n \in \te_\Psi(K) \cap \sigma \Gamma \sigma^{-1}$ such that for every
$\alpha \in \Pi \setminus \Psi$ we have: $\underset{n}{\lim}|\alpha(t_n^{-1})|_{{i}} = 0$
when $1 \leq i \leq s_1$, $\underset{n}{\lim}|\alpha(t_n)|_{i} = 0$ when $s_1 + 1 \leq i \leq s_2$, and $|\alpha(t_n)|_{i}$ is
bounded when $s_2+1 \leq i \leq r$.
\end{lem}

{\bf Proof.} The lemma follows from Proposition \ref{integers1'} and the commensurability of $\te(\OO)$ and $\te(K) \cap \sigma \Gamma \sigma^{-1}$.
\qed

\begin{lem}\label{apply units} Let $\Psi$, $\sigma$, $s_1$, and $s_2$ be as in the formulation of Lemma \ref{units}.
Also let $u_i \in \ve_\Psi^-(K_{i})$  if $1 \leq i \leq s_1$, $u_i \in \ve_\Psi(K_{i})$
if $s_1 + 1 \leq i \leq s_2$ and $g_i \in G_i$ if $s_2+1 \leq i \leq r$.
Then the closure of the orbit ${T_\Psi \pi(u_1 \sigma, \cdots, u_{s_2}\sigma, g_{s_2 + 1}, \cdots, g_r)}$ contains $\pi(\sigma, \cdots, \sigma, g_{s_2 + 1}, \cdots, g_r)$.
\end{lem}

{\bf Proof.} Let $t_n \in \te(K) \cap \sigma \Gamma \sigma^{-1}$ be as in Lemma \ref{units}. Passing to a subsequence we suppose
that for every $i > s_2 $ the projection of the sequence $t_n$ in $T_i$ is convergent. In view of the choice of $t_n$,
$t_n \pi(\sigma) = \pi(\sigma)$ and $\underset{n}{\lim} \ t_n u_i t_n^{-1} = e$ for all $1 \leq i \leq s_2$
(cf. (\ref{eq: uniprad of max parabolic}) and (\ref{eq: uniprad- of max parabolic})). Therefore
$$
\underset{n}{\lim} \ t_n\pi( u_1 \sigma, \cdots, u_{s_2}\sigma,  g_{s_2 + 1}, \cdots, g_r)
 = (e,\cdots, e ,h_{s_2 + 1} , \cdots, h_r)\pi(\sigma),
$$
where $h_i = \underset{n}{\lim} \  t_n g_i \sigma^{-1} t_n^{-1}, i > s_2$.
Using once again the convergence of $t_n$ in every $T_i$, $i > s_2 $, we get
$$
\underset{n}{\lim} \  t_n^{-1}(e,\cdots, e ,h_{s_2 + 1} , \cdots, h_r)\pi(\sigma) = \pi(\sigma, \cdots, \sigma, g_{s_2 + 1}, \cdots, g_r).
$$
\qed

\begin{lem}\label{lem3} Let $\Psi \subset \Pi$ and $g \in \G(K)$. Then
\begin{enumerate}
\item[(a)] $g = \omega z v_+ v_-$, where $\omega \in \mathcal{N}_{\G}(\te)(K)$,
$z \in \mathcal{Z}_{\G}(\te_\Psi)(K)$, $v_+ \in \ve_\Psi^+(K)$ and $v_- \in \ve_\Psi^-(K)$. Moreover,
$$
\mathcal{Z}_{\T_\Psi}(g) = \mathcal{Z}_{\T_\Psi}(v_-) \cap \mathcal{Z}_{\T_\Psi}(v_+) \cap \mathcal{Z}_{\T_\Psi}(\omega).
$$
\item[(b)] With $g = \omega z v_+ v_-$ as in $(a)$, suppose that $\dim \mathcal{Z}_{\T_\Psi}(g) \geq \dim \mathcal{Z}_{\T_\Psi}(\theta g)$
for every $\theta \in \mathcal{N}_{\G}(\te)$. Then
$$
\mathcal{Z}_{\T_\Psi}(\omega)^{\circ}\supset \mathcal{Z}_{\T_\Psi}(g)^{\circ}= \big(\mathcal{Z}_{\T_\Psi}(v_-)
\cap \mathcal{Z}_{\T_\Psi}(v_+)\big)^{\circ}.
$$
\end{enumerate}
\end{lem}
{\bf Proof.} It follows from the Bruhat decomposition  \cite[21.15]{Borel} and the structure of
the standard parabolic subgroups (see 2.3) that
\begin{eqnarray*}
& &
\G(K) = \ve_{\emptyset}^-(K)\mathcal{N}_{\G}(\te)(K)\ve_{\emptyset}^-(K) = \mathcal{N}_{\G}(\te)(K)\ve_{\emptyset}^+(K)\ve_{\Psi}^-(K) =
\\&& \mathcal{N}_{\G}(\te)(K)\ve_{\Psi}^+(K)\pe_{\Psi}^-(K) = \mathcal{N}_{\G}(\te)(K)\mathcal{Z}_{\G}(\te_\Psi)(K) \ve_{\Psi}^+(K)\ve_{\Psi}^-(K).
\end{eqnarray*}
So, $g = \omega z v_+ v_-$ as in the formulation of the lemma. Let $t \in \mathcal{Z}_{\T_\Psi}(g)$.
We have
$$
g = tgt^{-1} = \omega z v_+ v_- =  \omega (\omega^{-1}t\omega t^{-1} z)(tv_+t^{-1}) (tv_-t^{-1}).
$$
The product map $$\mathcal{Z}_{\G}(\te_\Psi)(K) \times \ve_{\Psi}^+(K) \times \ve_{\Psi}^-(K) \rightarrow \mathcal{Z}_{\G}(\te_\Psi)(K)
 \ve_{\Psi}^+(K) \ve_{\Psi}^-(K), (x,y,z) \mapsto xyz,$$ being bijective,  we obtain
 $$z = \omega^{-1}t\omega t^{-1} z,  \ v_+ = tv_+t^{-1} \ \textrm{and} \ v_- = tv_-t^{-1}.$$
We have proved that $\mathcal{Z}_{\T_\Psi}(g) \subset \mathcal{Z}_{\T_\Psi}(v_-) \cap \mathcal{Z}_{\T_\Psi}(v_+) \cap \mathcal{Z}_{\T_\Psi}(\omega)$.
Since the opposite inclusion is obvious, (a) is proved.

 The part (b) of the lemma follows immediately from (a). \qed

\begin{prop}\label{lem4} Let $g \in \G(K)$ be such that $\dim \mathcal{Z}_{\T}(g) \geq \dim \mathcal{Z}_{\T}(\theta g)$ for
all $\theta \in \mathcal{N}_{\G}(\te)$. Let $\Lambda$ be a subset of $\Phi$ and
$\mathbf{S} = \big(\underset{\alpha \in \Lambda}{\bigcap}\mathrm{ker} \ \alpha \big)^{\circ}$. Then there exist systems of
simple roots $\Pi$ and $\Pi'$ in $\Phi$ and subsets $\Psi \subset \Pi$ and $\Psi' \subset \Pi'$ with the following properties:
\begin{enumerate}
\item[(a)] $\mathbf{S} = \te_{\Psi} = \te_{\Psi'}$;
\item[(b)] $g  = \omega z v_+ v_-= \omega' z' v_-' v_+',$
where $\omega$ and $\omega' \in \mathcal{N}_{\G}(\te)(K)$, $z$ and $z' \in \mathcal{Z}_\G(\mathbf{S})(K)$, $v_+ \in \ve_\Psi(K)$, $v_- \in
\ve_\Psi^-(K)$, $v_+' \in \ve_{\Psi'}(K)$ and $v_-' \in \ve_{\Psi'}^-(K)$, and
$$
\mathcal{Z}_{\mathbf{S}}(g)^{\circ}= \mathcal{Z}_{\mathbf{S}}(v_+)^{\circ}= \mathcal{Z}_{\mathbf{S}}(v_+')^{\circ}.
$$
\end{enumerate}
\end{prop}
{\bf Proof.} Let us prove the existence of $\Psi$. (The proof of the existence of $\Psi'$ is analogous.)
 Fix $v \in \mathcal{S}$. We may (as we will) choose a $t \in \se(K_v)$ with the property: $|\alpha(t)|_v \neq 1$ for every root $\alpha$
which is not a
linear combination of roots from $\Lambda$. Applying Proposition \ref{horospherical}, fix a system of simple roots $\Pi$ and
a subset $\Psi$ of $\Pi$ such that
$\mathbf{S} = \te_{\Psi}$, $W^+(t) = \ve_\Psi(K_v)$, $W^-(t) = \ve_\Psi^-(K_v)$ and $\mathcal{Z}_\G(t) = \mathcal{Z}_\G(\te_{\Psi})$.
Let $g  = \omega z v_+ v_-$ as given by Lemma \ref{lem3}. Suppose that $t$ is chosen in such a way
that $\dim \mathcal{Z}_{\mathbf{S}}(v_+)$ is minimum possible. In view of Lemma \ref{lem3}(b), in order to show that $\Psi$ is as needed,
it is enough to prove that
$\mathcal{Z}_{\mathbf{S}}(v_+)^{\circ}\subset \mathcal{Z}_{\mathbf{S}}(v_-)^{\circ}$. Suppose by the contrary that
$\mathcal{Z}_{\mathbf{S}}(v_+)^{\circ}\nsubseteq \mathcal{Z}_{\mathbf{S}}(v_-)^{\circ}$. Pick a $t' \in \mathcal{Z}_{\mathbf{S}}(v_+)^{\circ}$
such that the subgroup generated by $t'$ is Zariski dense in $\mathcal{Z}_{\mathbf{S}}(v_+)^{\circ}$ and for every $K$-root $\beta$
either $\beta(t') = 1$ or $|\beta(t')|_v \neq 1$. Then $v_- = w_+ w_0 w_-$ where $w_+ \in W^+(t')$, $w_- \in W^-(t')$ and $w_0 \in \mathcal{Z}_\G (t')$.
(We use that $\ve_\Psi$ is directly spanned in any order by its subgroups $\mathbf{U}_{(\alpha)}$, see \cite[21.9]{Borel}.) Since $\mathcal{Z}_{\mathbf{S}}(v_+)^{\circ}\nsubseteq \mathcal{Z}_{\mathbf{S}}(v_-)^{\circ}$, we have that either $w_+ \neq e$ or $w_- \neq e$.
 Replacing, if necessary, $t'$ by $t'^{-1}$ we may (and will) suppose that $w_+ \neq e$. Put $\widetilde{t} = t t'^n$, $n \in \N$.
With $n$ chosen sufficiently large, $\widetilde{t}$ has the properties: $|\alpha(\widetilde{t})|_v \neq 1$ for every root $\alpha$ which is not a linear combination
of roots from $\Lambda$, $v_+w_+ \in W^+(\widetilde{t})$ and $w_0w_- \in W^-(\widetilde{t})$. In view of the choice of $w_+$,
$\mathcal{Z}_{\mathbf{S}}(v_+w_+) = \mathcal{Z}_{\mathbf{S}}(v_+) \cap \mathcal{Z}_{\mathbf{S}}(w_+)$. Since $w_+ \neq e$ and $t'$
centralizes $v_+$ but not $w_+$, we obtain
that $\dim \mathcal{Z}_{\mathbf{S}}(v_+w_+) < \dim \mathcal{Z}_{\mathbf{S}}(v_+)$ which contradicts the choice of $t$.
Therefore $\mathcal{Z}_{\mathbf{S}}(v_+)^{\circ}\subset \mathcal{Z}_{\mathbf{S}}(v_-)^{\circ}$ proving that $\Psi$ is as needed.
\qed

\begin{prop}\label{apply units+} Let $g = (g_1, \cdots, g_r) \in G$ where $g_i \in \G(K)$ for all $i$.
Let $\Psi \subset \Pi$ and $g_i g_r^{-1} = \omega_i z_i v_i^+ v_i^-$, where $\omega_i \in \mathcal{N}_{\G}(\te)(K)$,
$z_i \in \mathcal{Z}_{\G}(\te_\Psi)(K)$, $v_i^+ \in \ve_\Psi^+(K)$ and $v_i^- \in \ve_\Psi^-(K)$, for all $1 \leq i \leq r-1$.
Put $h_i = v_i^- \cdot g_r$, $1 \leq i < r$ and $h_r = g_r$. 
Then $\overline{T\pi(g)}$ contains the elements
$$\pi((g_1, \cdots, g_{i-1}, \omega_i z_i h_i, \cdots, \omega_{r-1} z_{r-1} h_{r-1}, h_r))$$
for all $1 \leq i \leq r-1$.
\end{prop}
{\bf Proof.} Assume, as we may, that all $\omega_i = 1$. In this case it is enough to prove that $\overline{T_\Psi \pi(g)}$,
 where $T_\Psi = \te_\Psi(K_{\mathcal{S}})$, contains
$\pi((g_1, \cdots, g_{i-1}, z_i h_i, \cdots, z_{r-1} h_{r-1}, h_r))$. Since $z_i$ centralize $T_\Psi$, without
restriction, we assume that all $z_i = e$. We will proceed by induction on $r-i$. Writing
$$
g_{r-1} = v_{r-1}^+(v_{r-1}^-g_r) \ \textrm{and} \ g_r = ({v_{r-1}^-})^{-1}(v_{r-1}^-g_r),
$$
 it follows from Lemma \ref{apply units} that
$\overline{T_\Psi \pi(g)}$ contains $\pi((g_1, \cdots, g_{r-2}, h_{r-1}, h_{r}))$. Suppose, by the induction hypothesis, that
$\pi((g_1, \cdots, g_i, h_{i+1}, \cdots, h_{r})) \in \overline{T_\Psi \pi(g)}$. Since $g_i = v_i^+ (v_i^-g_r)$ and
$h_{i+1} = v_{i+1}^-(v_i^-)^{-1}(v_i^-g_r)$, applying again
 Lemma \ref{apply units}, we obtain that $\pi((g_1, \cdots, g_{i-1}, h_{i}, \cdots, h_{r})) \in \overline{T_\Psi \pi(g)}$.
\qed

\subsection{Definition of $h_1{H}_1\pi(e)$ in (\ref{skeezing}).}
Let
$g = (g_1, \cdots, g_r) \in G$ and $T\pi(g)$ be a locally divergent orbit.
According to Theorem \ref{ldo},  $g = zg'$ where $z \in \mathcal{Z}_{G}(T)$, $g' = (g'_1, \cdots, g'_r) \in G$ and all $g'_i \in \G(K)$.
Clearly, $zT\pi(g') = T\pi(g)$ where the orbit $T\pi(g')$ is locally divergent. Hence we may (and will) assume that all $g_i \in \G(K)$.

Next choose $\omega_i \in \mathcal{N}_{\G}(\te)(K), 1 \leq i \leq r-1,$ in such a way that
$\dim \overset{r-1}{\underset{i=1}{\bigcap}} \mathcal{Z}_{\te}(\omega_i g_i g_r^{-1})$ is maximum possible. Let
$\mathbf{H}'_1 = \mathcal{Z}_{\G}\big((\overset{r-1}{\underset{i=1}{\bigcap}} \mathcal{Z}_{\te}(\omega_i g_i g_r^{-1}))^{\circ}\big)$
and $H_1' = \mathbf{H}'_1(K_\mathcal{S})$.
In view of \cite[11.12]{Borel}, $\mathbf{H}'_1$ is a connected reductive $K$-group containing $\te$
and $(\overset{r-1}{\underset{i=1}{\bigcap}} \mathcal{Z}_{\te}(\omega_i g_i g_r^{-1}))^{\circ}$ is its connected center.
We put
$$
\mathbf{H}_1 = g_r^{-1}\mathbf{H}'_1 g_r.
$$
Let us show that
$$T\pi(g) \subset h_1 H_1 \pi(e),$$
where $h_1 = (\omega_1^{-1} g_r, \cdots , \omega_{r-1}^{-1}g_{r}, g_r)$ and $H_1 = \mathbf{H}_1(K_{\mathcal{S}})$.
Indeed, since $\omega_ig_ig_r^{-1} \in \mathbf{H}_1'$ for all $i$, $(\omega_1, \cdots, \omega_{r-1}, e)$ normalizes $T$ and $T \subset {H}_1'$, we have
\begin{eqnarray*}
& & T\pi(g) \subset (\omega_1, \cdots, \omega_{r-1}, e)^{-1}H_1'(\omega_1, \cdots, \omega_{r-1}, e)\pi(g) =
\\ & & (\omega_1, \cdots, \omega_{r-1}, e)^{-1}H_1'(\omega_1 g_1 g_r^{-1}, \cdots, \omega_{r-1} g_{r-1} g_r^{-1}, e)\pi(g_r)= \\ & &
(\omega_1, \cdots, \omega_{r-1}, e)^{-1}H_1'\pi(g_r) =  h_1 H_1 \pi(e).
\end{eqnarray*}
Remark that the orbit $H_1'\pi(g_r)$ is closed and $T$-invariant.
Therefore
$h_1 H_1 \pi(e)$ is closed and $T$-invariant as in the formulation of Theorem \ref{r>2}.
\subsection{Reducing the proof of Theorem \ref{r>2}.} The proof of the existence of $h_2{H}_2\pi(e)$ as in (\ref{skeezing})
represents the main part of the proof of  Theorem \ref{r>2}. 

With the notation of 6.2, we have
\begin{equation}
\label{def H_1}
\mathbf{H}'_1 = \mathcal{Z}_{\G}\big((\overset{r-1}{\underset{i=1}{\bigcap}} \mathcal{Z}_{\te}(\omega'_i\omega_i g_i g_r^{-1}))^{\circ}\big)
\end{equation}
for all choices of $\omega'_i \in \mathcal{N}_{\mathbf{H}'_1}(\te)$. 

Note that $g_r^{-1}\te g_r  \subset \mathbf{H}_1$, $g_r^{-1}\omega_i g_i \in  \mathbf{H}_1$ for all
$i$, $\mathcal{Z}(\mathbf{H}_1)^\circ = (\overset{r-1}{\underset{i=1}{\bigcap}} \mathcal{Z}_{\te}(g_r^{-1}\omega_i g_i))^\circ$,
$$h_1^{-1}T\pi(g) = (g_r^{-1} T g_r)\pi((g_r^{-1} \omega_1 g_1, \cdots, g_r^{-1} \omega_{r-1}g_{r-1}, e)) \subset  H_1 \pi(e)$$
and $\mathbf{H}_1$ is an almost direct product over $K$ of $\mathcal{Z}(\mathbf{H}_1)$ and $\mathcal{D}(\mathbf{H}_1)$.
Therefore the locally divergent orbit $h_1^{-1}T\pi(g)$ gives rise to a locally divergent orbit on the quotient of
$\mathcal{D}(\mathbf{H}_1)(K_{\mathcal{S}})$ by an arithmetic subgroup reducing the proof of Theorem \ref{r>2} to the following case:

$(\ast)$  \textit{all $g_i \in \G(K)$, $g_r = e$, $\underset{i}{\bigcap} \mathcal{Z}_{\te}(\omega_i g_i)$
is finite for all choices of $\omega_i \in \mathcal{N}_{\G}(\te)(K)$ and, in view of (\ref{def H_1}),
$\dim \mathcal{Z}_{\T}(g_i) \geq \dim \mathcal{Z}_{\T}(\omega g_i)$ whenever $1 \leq i < r$ and $\omega \in \mathcal{N}_{\G}(\te)(K)$.}

It remains to prove that, under the conditions of $(\ast)$, there exists a semisimple $K$-subgroup $\mathbf{H}$ of $\G$ with
$\mathrm{rank}_K (\G) = \mathrm{rank}_K (\mathbf{H})$, a subgroup of finite index $H$ in $\mathbf{H}(K_{\mathcal{S}})$, and $h \in \G(K)$ such that
 $h H \pi(e)$ is $T$-invariant and
$$
h H \pi(e) \subset \overline{T\pi(g)}.
$$

\subsection{Special elements in $\overline{T\pi(g)}$.} We will suppose up to the end of
section 6.8 that the conditions of $(\ast)$ are fulfilled.

\begin{prop}\label{prop4} For every $j$, $1 \leq j \leq r$, $\overline{T\pi(g)}$ contains an element of the form
$\omega(e, \cdots, e, \underset{j}{u},e, \cdots, e)\pi(h)$, where $h \in \G(K)$, $\omega \in \mathcal{N}_G(T)$,
and $u$ is a unipotent element in $\G(K)$ such that
 $\mathcal{Z}_{\T}(u)$ is finite.
\end{prop}

{\bf Proof.}
First, we will prove the proposition in the particular case when $\mathcal{Z}_\T(g_i)$ is finite for some $i$. Let $i = 1$.
By Proposition \ref{lem4} there exists a system of simple roots $\Pi$ such that every
 $g_i = z_i u_i^+u_i^-$, where $u_i^+ \in \ve_{\emptyset}(K)$,  $u_i^- \in \ve_{\emptyset}^-(K)$,
and $z_i \in \mathcal{N}_\G(\T)(K)$, and, moreover, $\mathcal{Z}_\te(u_1^+)$ is finite. Shifting $g$ from the left by an appropriate
element from $\mathcal{N}_G(T)$ we may assume that all $z_i = e$. It follows from Proposition \ref{apply units+} that $\overline{T\pi(g)}$
contains
$$\pi((u_1^{+}u^{-}_{1}, u^{-}_{2}, \cdots, u^{-}_{2})) = (u_1^{+}u^{-}_{1}(u^{-}_{2})^{-1}, e, \cdots, e)\pi(u_{2}^-).$$
Put $h = u_{2}^-$. By Proposition \ref{lem4}
 there exist opposite minimal parabolic $K$-subgroups  $\widetilde{\pe}^-$ and $\widetilde{\pe}^-$
 containing $\te$ such that $u_1^{+}u^{-}_{1}(u^{-}_{2})^{-1} = zu^-u$, where $z \in \mathcal{N}_{\G}(\T)(K)$,
$ u^- \in \mathcal{R}_u(\widetilde{\pe}^-)(K)$, $u \in \mathcal{R}_u(\widetilde{\pe}^+)(K)$ and
$ \mathcal{Z}_{\T}(u)$ is finite. We may suppose that $z = e$. Let $1 < j \leq r$. Writing
$$
(u^-u,  e, \cdots, e)\pi(h) = (u^-u,  e, \cdots, e, \underset{j}{(u)^{-1}u}, e, \cdots, e)\pi(h),
$$
Lemma \ref{apply units} implies that $\overline{T\pi(g)}$ contains $(u,  e, \cdots, \underset{j}{u}, \cdots, e)\pi(h)$.
Using the assumption $r > 2$ and applying once again Lemma \ref{apply units}, we obtain that $\overline{T\pi(g)}$ contains both
$(e,\cdots, \underset{j}{u}, \cdots, e)\pi(h)$ and $(u, e, \cdots, e)\pi(h)$.

In order to reduce the proof of the proposition to the case considered above,
it is enough to prove that if $i \neq j$, say $i =1$ and $j = 2$, then
$\overline{T\pi(g)}$ contains an element $\pi((g_1', g_2',g_3,$ $\cdots, g_r))$ such that $g_1'$ and $g_2' \in \G(K)$,
$\dim \mathcal{Z}_\T(g'_i) \geq \dim \mathcal{Z}_\T(\omega g'_i)$ for all $\omega \in \mathcal{N}_\G(\T)$ and $i \in \{1,2\}$, and
$$
\mathcal{Z}_{\T}(g_1')^\circ = \big(\mathcal{Z}_{\T}(g_1) \cap \mathcal{Z}_{\T}(g_2)\big)^\circ.
$$
We choose $g_1' = g_1$ if $\mathcal{Z}_{\T}(g_1)^\circ \subset \mathcal{Z}_{\T}(g_2)$. Suppose that
$\mathcal{Z}_{\T}(g_1)^\circ  \nsubseteq \mathcal{Z}_{\T}(g_2)$. By Lemma \ref{lem3} and Proposition \ref{lem4} there exists a
system of simple roots $\Pi$ and $\Psi \subset \Pi$ such that
$\te_\Psi = \mathcal{Z}_{\T}(g_1)^\circ$, $g_2  = \omega zv_- v_+$, where $\omega \in \mathcal{N}_{\G}(\mathbf{T})(K)$,
$z \in \mathcal{Z}_\G(\mathbf{T}_\Psi)(K)$,
$v_+ \in \ve_\Psi(K)$ and $v_- \in \ve_\Psi^-(K)$, and
$$
\mathcal{Z}_{\mathbf{T}_\Psi}(g_2)^\circ = \mathcal{Z}_{\mathbf{T}_\Psi}(v_+)^\circ.
$$
Since $\pi(g) = (g_1(v_+)^{-1}v_+, \omega z v_- v_+, g_3, \cdots, g_r)\pi(e)$,
Lemma \ref{apply units} implies that $\overline{T\pi(g)}$ contains $(g_1v_+, \omega z v_+, g_3, \cdots, g_r)\pi(e)$.
It is clear that
 \begin{eqnarray*}
&& \mathcal{Z}_{\mathbf{T}}(g_1v_+)^\circ = (\mathcal{Z}_{\mathbf{T}}(g_1) \cap \mathcal{Z}_{\mathbf{T}}(v_+))^\circ = \\ &&
 (\mathcal{Z}_{\mathbf{T}}(g_1) \cap \mathcal{Z}_{\mathbf{T}_\Psi}(g_2))^\circ =(\mathcal{Z}_{\mathbf{T}}(g_1) \cap \mathcal{Z}_{\mathbf{T}}(g_2))^\circ,
 \end{eqnarray*}
compleating the proof.
\qed

\medskip

Proposition \ref{prop4} is strengthen as follows.

\begin{cor}\label{abelian} With the notation and assumptions of Proposition \ref{prop4}, $\overline{T\pi(g)}$ contains an element of the form
$\omega(u, e, \cdots, e)\pi(h)$, where $\omega \in \mathcal{N}_G(T)$, $h \in \G(K)$, $u$ belongs to an
abelian unipotent subgroup of $\G(K)$ normalized by $\mathbf{T}(K)$, and
$\mathcal{Z}_{\T}(u)$ is finite.
\end{cor}

We need the following.

\begin{lem}\label{convex} Consider the $\Q$-vector space $\Q^n$ endowed with the standard scalar product:
$((x_1, \cdots, x_n),(y_1, \cdots, y_n)) \bydefn \underset{i}{\sum}x_i y_i$. Let $v_1, \cdots, v_m$ be pairwise
non-proportional vectors in $\Q^n$  and $v \in \Q^n$ be such that $(v_i, v) > 0$ for all $1 \leq i \leq m$.
Let $\mathcal{C} = \{\overset{m}{\underset{i=1}{\sum}}a_i v_i | a_i \in \Q, a_i \geq 0\}$. Suppose that
$m > n$ and the interior of the cone $\mathcal{C}$ with respect to the topology on $\Q^n$ induced by $(\cdot,\cdot)$
is not empty. Then there exist $1 \leq {i_\circ} \leq m$ and $w \in \Q^n$ such that $\{v_i | i \neq i_\circ\}$ contains a basis of $\Q^n$,
$(w,v_{i_\circ}) < 0$, and
$(w,v_{i}) > 0$ for all $i \neq i_\circ$.
\end{lem}
{\bf Proof.} Put $\mathcal{C}_{\R} = \{\overset{m}{\underset{i=1}{\sum}}a_i v_i | a_i \in \R, a_i \geq 0\}$.
Rearranging the indices of $v_i$, we may assume that $\{v_1, \cdots, v_{m_1}\}$ is a minimal subset of $\{v_1, \cdots, v_{m}\}$
such that $\mathcal{C}_{\R} = \{\overset{m_1}{\underset{i=1}{\sum}}a_i v_i | a_i \in \R, a_i \geq 0\}$ and, moreover, $\{v_1, \cdots, v_{n}\}$
is a bases of $\R^n$.  If $v_{n+1} = \overset{n}{\underset{i=1}{\sum}}b_i v_i$ then at least one of the $b_i$ is strictly positive. Let $b_1 > 0$.  Let
$\mathcal{C}'_{\R} = \{\overset{m}{\underset{i=2}{\sum}}a_i v_i | a_i \in \R, a_i \geq 0\}$. Then $\mathcal{C}'_{\R} \varsubsetneq \mathcal{C}_{\R}$,
$v_1 \in \mathcal{C}_{\R} \setminus \mathcal{C}'_{\R}$ and $\{v_2, \cdots, v_{m}\}$ contains a basis of $\Q^n$. Using, for example, the Hahn-Banach theorem about separation of convex subsets of an
 affine space by hyperplans \cite[11.4.1]{Berget}, one can prove by a standard argument the existence of $w \in \R^n$ such that $(w,v_{1}) < 0$ and
$(w,v_{i}) > 0$ for all $i > 1$. Since $\Q^n$ is dense in $\R^n$, we can choose $w$ in $\Q^n$. \qed

\medskip

{\bf Proof of Corollary \ref{abelian}.} By Proposition \ref{prop4} $\overline{T\pi(g)}$ contains an element
$\omega(u, e, \cdots, e)\pi(h)$, where $\omega \in \mathcal{N}_G(T)$, $h \in \G(K)$ and $u \in \G(K)$ is a
 unipotent element such that $\mathcal{Z}_{\T}(u)$ is finite. Let $\ve$ be the minimal $\T$-invariant unipotent $K$-subgroup of $\G$ containing $u$.
 We assume, with no loss of generality, that the element $u$ with the above properties is such that $\dim \V$ is minimal possible.
 It remains to prove that $\ve$ is abelian.
 The proof is easily reduced to the case when $\omega = e$. Suppose by the contrary that $\ve$ is not abelian.
 There exists a system of positive roots $\Phi^+$ in $\Phi$
such that the corresponding to $\Phi^+$ maximal unipotent $K$-subgroup of $\G$ contains $\ve$. Let $\Phi^+_{nd}$ be the subset of non-divisible roots
in $\Phi^+$ and $\{\alpha_1, \cdots, \alpha_m\} = \{\alpha \in \Phi^+_{nd}: \U_{(\alpha)} \cap \ve \neq \{e\}\}$ where $\U_{(\alpha)}$ is
the corresponding to $\alpha$ subgroup (see 2.3). Put $\ve_{(\alpha_i)} = \U_{(\alpha_i)} \cap \ve$. Rearranging $\{\alpha_1, \cdots, \alpha_m\}$,
we assume that $l \leq m$ is such that every $\alpha_i$, $i > l$, is a linear combination with strictly positive coefficients
of at least two different roots from $\{\alpha_1, \cdots, \alpha_l\}$ and no one of the roots in $\{\alpha_1, \cdots, \alpha_l\}$ has this property.
It follows from the standard rule $[\mathfrak{g}_\alpha, \mathfrak{g}_\alpha] \subset \mathfrak{g}_{\alpha+\beta}$ that
the product in any order of all $\ve_{(\alpha_i)}$, $i > l$, is a normal subgroup of $\ve$ which we denote by $\ve'$. Also, for every $i \leq l$,
 the product in any order of all
$\ve_{(\alpha_j)}$, where $1 \leq j \leq m$ and $j \neq i$, is a normal subgroup of $\ve$ which we denote by $\ve_i'$. Since $\ve /\ve_i'$ is isomorphic to
$\ve_{(\alpha_i)}$ and the group $\U_{(\alpha)}$ is abelian if ${(\alpha)} = \{\alpha\}$ or metabelian with center $\U_{2\alpha}$ if
${(\alpha)} = \{\alpha, 2\alpha\}$ \cite[21.10]{Borel}, it follows from the definition of $\ve$ that $\ve_{(\alpha_i)}$ is abelian if $i \leq l$. Hence
$\ve'$ contains the derived subgroup $\mathfrak{D}(\ve)$ of $\ve$.

Let $u = u_1 \cdots u_m$ where  $u_i \in \ve_{(\alpha_i)}(K)$. Suppose on the contrary that $\mathfrak{D}(\ve)$ is not trivial.
Then $u_i \neq e$ for all $1 \leq i \leq l$ and $u_j \neq e$ for some $j > l$. Consider the $\Q$-vector space $X(\T)\underset{\Z}{\bigotimes} \Q$ and
the cone $\mathcal{C} = \{\overset{m}{\underset{i=1}{\sum}}a_i \alpha_i | a_i \in \Q, a_i \geq 0\}$. Using Lemma \ref{convex} and the natural pairing between
the group of characters of $\te$ and the group of the multiplicative one-parameter subgroups in $\te$ \cite[8.6]{Borel}, we find
$1 \leq i \leq l$, say $i = 1$, and $t \in \te_1$ such that $\underset{n \rightarrow +\infty}{\lim} t^n u_1 t^{-n} = e$ in $G_1$ and
$\underset{n \rightarrow -\infty}{\lim} t^n u_i t^{-n} = e$ in $G_1$ for all $i > 1$. Put $u' = u_2 \cdots u_m$. It follows from
Proposition \ref{apply units+} that $\overline{T\pi(g)}$ contains $(u', e, \cdots, e)\pi(h)$. It remains to note that $\mathcal{Z}_{\te}(u')$ is finite
and $u'$ is contained in a proper $\T$-invariant $K$-subgroup of $\ve$ which is a contradiction.\qed

\subsection{Unipotent orbits on $\overline{T\pi(g)}$.} Further on, we denote by $\mathcal{O}'$ a subring of finite index in $\mathcal{O}$ and put
$\mathcal{O}'_\infty = \mathcal{O}_\infty \cap \mathcal{O}'$.

Some propositions in $\mathcal{S}$-adic setting
will be deduced from their archimedean analogs when $\mathcal{S} = \mathcal{S}_\infty$. For this purpose we need

\begin{lem}\label{approximation} Let $\ve$ be a unipotent $K$-algebraic group and $\U$ be its $K$-subgroup. Put $U = \U(K_\mathcal{S})$
and $U_\infty = \U(K_\infty)$.  Let $M$ be a subset of
$U_{\infty}$ such that $\overline{M \ve(\mathcal{O}'_{\infty})} = U_{\infty}\ve(\mathcal{O}'_{\infty})$. Then
$$
\overline{M \ve(\mathcal{O}')} = U\ve(\mathcal{O}').
$$
\end{lem}
{\bf Proof.} By the strong approximation for unipotent groups (see, for example, \cite[\S 7.1, Corollary]{PR}), we have
that $U = \overline{U_{\infty}\U(\mathcal{O}')}$. Using that $\U(\OO') \subset \ve(\OO')$, $\ve(\OO'_\infty) \subset \ve(\OO')$ and $U\ve(\OO')$ is closed,
we get
$$
U\ve(\OO') = \overline{U_{\infty}\ve(\mathcal{O}')} =
\overline{U_{\infty}\ve(\mathcal{O}'_\infty)\ve(\mathcal{O}')}=\overline{M\ve(\mathcal{O}'_\infty)\ve(\mathcal{O}')}=\overline{M\ve(\OO')}.
$$
\qed

\begin{prop}\label{prop5} We suppose that
$K$ is not a $\mathrm{CM}$-field and the completion $K_1$ is archimedean.
Let $\mathbf{V}$ be an abelian unipotent $K$-subgroup of $\G$ normalized by $\mathbf{T}$.
Let $u \in \mathbf{V}(K)$ and
 $\mathcal{Z}_{\T}(u)$ be finite.
Then there exists a $K$-subgroup $\mathbf{U}$ of $\mathbf{V}$ which is $\mathbf{T}$-invariant, contains $u$, and
\begin{equation}
\label{abelian1}
U\pi(e) = \overline{\mathbf{T}(\mathcal{O}')(u, e, \cdots, e)\pi(e)} = \overline{\{(tut^{-1},e, \cdots,e)\pi(e) : t \in \mathbf{T}(\mathcal{O}')\}},
\end{equation}
where $U = \mathbf{U}(K_{\mathcal{S}})$.
\end{prop}

We denote by $L$ the closure of the projection of $\mathcal{O}'^*$ in $K_1^*$. We will consider separately the
two cases: $L$ has finite index in  $K_1^*$ and $L$ has infinite index in $K_1^*$. (Clearly, $L$ has finite index in  $K_1^*$
if and only if either $L = K_1^* = \C^*$ or $K_1 = \R$ and $L$ contains the strictly positive reals.)

{\bf Proof of Proposition \ref{prop5} when $L$ has finite index in  $K_1^*$.} With the notation of 2.3, there exists an order $\Phi^+$ of the set of $K$-roots
with respect to $\te$ such that $\ve \subset \ve_\emptyset$. Identifying $T_1$ with $(K_1^*)^{\dim \te}$,
the map $T_1 \rightarrow \ve(K_1), t \mapsto tut^{-1},$ can be regarded as the restriction to $(K_1^*)^{\dim \te}$ of a
polynomial map $K_1^{\dim \te} \rightarrow \ve(K_1)$.
Let $\pi_\infty : V_\infty \mapsto V_\infty/\ve(\mathcal{O'_\infty})$ where $V_\infty = \ve(K_\infty)$, be the natural projection.
Remark that $V_\infty/\ve(\mathcal{O'_\infty})$ is a usual topological compact torus.
By the polynomial orbit rigidity for tori (see \cite[Theorem 8]{Weyl} or \cite[Corollary 1.2]{Shah} for a general recent result),
there exists a $\te$-invariant $K$-subgroup $\mathbf{U}$ of $\mathbf{V}$ such that
$$
U_\infty \pi_\infty(e) = \overline{\{(tut^{-1},e, \cdots,e)\pi_\infty(e) : t \in T_1\}},
$$
where $U_\infty = \mathbf{U}(K_\infty)$. Since $L = K_1^*$, (\ref{abelian1}) follows from Lemma \ref{approximation}
applied to $M \bydefn \{(tut^{-1},e, \cdots,e) : t \in \te(\mathcal{O}')\} \subset U_\infty$. \qed

\medskip

{\bf Proof of Proposition \ref{prop5} when $L$ has infinite index in $K_1^*$.} Since $K$ is not a $\mathrm{CM}$-field,
Proposition \ref{integers1} implies that $K_1 = \C$, $\dim L = 1$, $L \neq \R_+$, and $K_1^*/L$ is compact.
Up to a subgroup of finite index there are two possibilities for $L$: (a) $L$ is a direct product of the unit circle group
$S^1$ and an infinite cyclic group, i.e. $L = \{e^{2 \pi n\alpha  + \mathfrak{i} t}: n \in \Z, 0 \leq t < 2\pi\}$,
where $\mathfrak{i}^2 = -1$ and $\alpha \in \R^*$, and (b) $L$ is a spiral, i.e.
$L = \{e^{(\alpha + \mathfrak{i})(t +2\pi n)}: n \in \Z, 0 \leq t < 2\pi\}$ where $\alpha \in \R^*$.
Further on, we denote by $L$ a subgroup of $K_1^*$ satisfying (a) or (b).
With the above notation, the case $\alpha < 0$ being analogous to the
case $\alpha > 0$, we will suppose that $\alpha > 0$. In order to treat (a) and (b) simultaneously,
we write
 $L = \{e^{2 \pi n \alpha+(\widetilde{\alpha} + \mathfrak{i})t}: n \in \Z, 0 \leq t < 2 \pi\}$
where $\widetilde{\alpha} = 0$ in case (a) and $\widetilde{\alpha} = \alpha$ in case (b).
(We use the equality $e^{2 \pi n \alpha+({\alpha} + \mathfrak{i})t} = e^{(\alpha + \mathfrak{i})(t +2\pi n)}$.)

Fix an ordering $\Phi^+$ of the root system $\Phi$ such that $\ve \subset \ve_\emptyset$. The group $\ve$ is $K$-isomorphic to,
and will be identified with, a $K$-vector space. There exist pairwise different positive roots $\alpha_1, \cdots, \alpha_l$ such that
$\ve = \overset{l}{\underset{i}{\bigoplus}} \ve_{i}$ where $\ve_{i}$ is a weight subspace with weight $\alpha_i$ for the action of $\te$
on $\ve$. Since $\alpha_1, \cdots, \alpha_l$ are pairwise different and positive there exists a $1$-parameter group $\lambda: \mathbf{GL}_1 \rightarrow \te$
such that $\alpha_i \circ \lambda (t) = t^{n_i}$ where $n_i$ are pairwise different positive integers. With
$u$ as in the formulation of the proposition let $u = \underset{i}{\sum}u_i$ where $u_i \in \ve_{i}(K)$. Denote by $\mathbf{U}$ the subspace of $\ve$
spanned by all $u_i$. Let $U = \mathbf{U}(K_\mathcal{S})$. Then $U\pi(e)$ is closed containing ${\mathbf{T}(\mathcal{O}')(u, e, \cdots, e)\pi(e)}$. The
opposite inclusion (and the proof of Proposition \ref{prop5}) follows from the next lemma.

\begin{lem}
\label{Wey2} Let $\mathbf{GL}_1$ act $K$-rationally on a finite dimensional $K$-vector space
$\mathbf{U}$ and $\mathbf{U} = \overset{l}{\underset{i}{\bigoplus}} \mathbf{U}_{\lambda_i}$
be the decomposition of $\mathbf{U}$ as a sum of weight sub-spaces $\mathbf{U}_{\lambda_i}$ with weights
$\lambda_i(t) = t^{n_i}$. Suppose that $r > 2$, $K$ is not a $\mathrm{CM}$-field, $n_i$ are
pairwise different positive integers, and every $\mathbf{U}_{\lambda_i}$ is spanned by an $u_i \in \mathbf{U}_{\lambda_i}(K) \setminus \{0\}$.
Then for every real $C > 1$, we have
$${U} = \overline{\{(\underset{i}{\sum}\lambda_i(a)u_i,0, \cdots,0): a \in L, |a|_1 \geq C\}+ \mathbf{U}(\mathcal{O}')},$$
where $U = \mathbf{U}(K_{\mathcal{S}})$.
\end{lem}

{\bf Proof.} In the course of the proof, given $\theta \in [0, 2\pi)$ and $b < c$, we denote $\R_\theta = \{r e^{\mathfrak{i}\theta}: r \in \R \}$ and
$[b,c]_\theta = \{r e^{\mathfrak{i}\theta}: a < r < b \}$. Both $\R_\theta$ and $[b,c]_\theta$ are imbedded in $K_1$ and the latter is one of the factors
in the direct product $K_\mathcal{S}$ (or $K_\infty$).

If $\sigma \in \mathrm{GL}(\mathbf{U}(K))$ then $\sigma(\mathbf{U}(\OO'))$ is commensurable with $\mathbf{U}(\OO')$. Since
there is no restriction on the choice of the subring of finite index $\OO'$ in $\OO$, we may (and will) identify $\mathbf{U}(K)$ with $K^l$.
and $u_i$ with the standard basis of $K^l$. 
The projection of $K_1$ into $K_{\mathcal{S}}/\mathcal{O}$ being dense, it follows from \cite[Corollary 1.2]{Shah} (or \cite{Weyl})
that for every $C > 1$
$$K_{\mathcal{S}}^l = \overline{\{(\underset{i}{\sum}\lambda_i(a)u_i,0, \cdots,0): a \in K_1, |a|_1 \geq C\}+ \mathcal{O}'^l}.$$
In view of Lemma \ref{approximation}, it is enough to prove that
$$U_\infty  = \overline{\{(\underset{i}{\sum}\lambda_i(a)u_i,0, \cdots,0): a \in L, |a|_1 \geq C\}+\mathbf{U}(\mathcal{O}'_\infty)},$$
where $U_{\infty} = \mathbf{U}(K_\infty)$.

It is enough to consider the case when $0 < n_1 < n_2 < \cdots < n_l$.
For every $i$ we introduce a parametric curve $f_i: [0, 2\pi) \rightarrow K^*_1$, $t \mapsto e^{(\widetilde{\alpha} + \mathfrak{i})n_i t}$.
Recall that $\mathcal{O}'_\infty$ is a group of finite type, the diagonal imbedding of $\mathcal{O}'_\infty$ in $K_\infty$ is a lattice and
$\overline{K_1 + \mathcal{O}'_\infty} = K_\infty$. Therefore $\overline{\R_\theta + \mathcal{O}'_\infty}$, $0 \leq \theta < 2\pi$,
is a subspace of the real vector space $K_\infty$ and the set of all $\theta \in [0, 2\pi)$ such that
$\overline{\R_\theta + \mathcal{O}'_\infty} \varsubsetneq K_\infty$ is countable.
Note that the tangent line at $t$ of each
curve $f_i$ runs over all directions when $0 \leq t < 2\pi$. By the above there exists $0 \leq \psi < 2\pi$ such that
if the tangent line at $\psi$ of the curve $f_i$ is parallel to ${\R_{\theta_i}}$, $0 \leq \theta_i < 2\pi$,
then $\overline{\R_{\theta_i} + \mathcal{O}'_\infty} = K_\infty$ for all $1 \leq i \leq l$.

For every $n \in \mathbb{N}_+$, let $$F_n: [0, 2\pi) \rightarrow K_{\infty}^l, \ \ t \mapsto ((e^{2\pi n n_1 \alpha}f_1(t), \cdots, 0), \cdots, (e^{2\pi n n_l \alpha}f_l(t), \cdots, 0)).$$

Let $\varepsilon > 0$. A subset $M$ of $K_{\infty}^l/{\mathcal{O}'}_{\infty}^l$ will be called {\it $\varepsilon$-dense} if
the $\varepsilon$-neighborhood of any point in $K_{\infty}^l/{\mathcal{O}'}_{\infty}^l$ contains an element from $M$.
(As usual, $K_{\infty}^l/{\mathcal{O}'}_{\infty}^l$ is endowed with a metrics induced by the standard metrics on
$K_{\infty}$ considered as a real vector space.) Remark that since $K_{\infty}^l/{\mathcal{O}'}_{\infty}^l$ is a compact torus if $M$ is $\varepsilon$-dense
then every shift of $M$ by an element of $K_{\infty}^l/{\mathcal{O}'}_{\infty}^l$ is also $\varepsilon$-dense.

Now the lemma follows from the next

\medskip

{\bf Claim.} {\it With $\psi$ and $F_n$ as above, let $\varepsilon > 0$. There exist reals $A_\varepsilon > 0$ and $b_\varepsilon > 0$ such that if
$\psi - b_\varepsilon \leq c < d \leq \psi + b_\varepsilon$ and $e^{2\pi n n_1}(d-c) > A_\varepsilon$, where $c$ and $d \in \R$ and $n \in \N_+$, then
$$
\{F_n(t) + {\mathcal{O}'}_{\infty}^l | c \leq t \leq d\}
$$
is $\varepsilon$-dense in $K_{\infty}^l/{\mathcal{O}'}_{\infty}^l$.}

\medskip

We will prove the claim by induction on $l$.
Let $l = 1$, i.e. $F_n: [0, 2\pi) \rightarrow K_{\infty}, \ \ t \mapsto (e^{2\pi n n_1 \alpha}f_1(t), 0, \cdots, 0)$, where
$f_1(t) = e^{(\widetilde{\alpha} + \mathfrak{i})n_1 t}$.
It follows from the choice of $\psi$ that there exists a real $B_\varepsilon > 0$ such that the projection of
$[0,B_\varepsilon]_{\theta_1}$ into $K_\infty/\mathcal{O}'_\infty$ is $\frac{\varepsilon}{2}$-dense. Every shift of
$[0,B_\varepsilon]_{\theta_1} + \mathcal{O}_\infty$ by an element from $K_\infty$ is also $\frac{\varepsilon}{2}$-dense in
$K_\infty/\mathcal{O}'_\infty$. Choosing $A_\varepsilon > 0$ sufficiently large and $b_\varepsilon > 0$ sufficiently close to $0$ we get that if
$n$ is such that $e^{2\pi n n_1}(2b_\varepsilon) > A_\varepsilon$ then the length of the curve
$\{F_{n}(t)| \psi - b_\varepsilon \leq t \leq \psi + b_\varepsilon\}$
is greater than $B_\varepsilon$ and if $I$ is any connected piece of this curve of length $B_\varepsilon$ then $I$ is $\frac{\varepsilon}{2}$-close
(with respect to the Hausdorff metrics on $\C$) to a shift of $[0,B_\varepsilon]_{\theta_1}$. Hence the projection of $I$ into $K_\infty/\mathcal{O}'_\infty$
is $\varepsilon$-dense completing the proof for $l = 1$.

Suppose that $l > 1$ and the claim is valid for  $l - 1$. Let
$$\widetilde{F}_n(t) = ((e^{2\pi n n_1 \alpha}f_1(t), \cdots, 0), \cdots, (e^{2\pi n n_{l-1} \alpha}f_{l-1}(t), \cdots, 0)).$$
By the induction hypothesis for  $l - 1$ and the validity of the Claim for $K_\infty/\mathcal{O}'_\infty$, given $\varepsilon > 0$
there exist positive reals $A_\varepsilon$ and $b_\varepsilon$ such that if $\psi -b_\varepsilon \leq c < d \leq \psi + b_\varepsilon$
and $e^{2\pi n n_1}(d-c) \geq A_\varepsilon$ for some $n \in \N_+$ then
$\{\widetilde{F}_n(t) + {\mathcal{O}'}_{\infty}^{l-1}| c \leq t \leq d\}$
is $\varepsilon$-dense in $K_{\infty}^{l-1}/{\mathcal{O}'}_{\infty}^{l-1}$ and
$$
\{(e^{2\pi n n_l \alpha}f_l(t),0 \cdots, 0) + \mathcal{O}'_\infty| c' \leq t \leq d'\}
$$
is $\varepsilon$-dense in $K_\infty/\mathcal{O}_\infty$ whenever $c < c' < d' < d$ and $e^{2\pi n n_l}(d'-c') \geq A_\varepsilon$.

Further on, given $c_* < d_*$, we define the length of the parametric curve $\{\widetilde{F}_n(t)| c_* \leq t \leq d_*\} \subset K_{\infty}^{l-1}$
as the maximum of the lengths of the curves $\{e^{2\pi n n_i \alpha}f_i(t) | c_* \leq t \leq d_*\} \subset \C$, $1 \leq i \leq l-1$.
With $b_\varepsilon$, $A_\varepsilon$ and $n$ as above, let $x \in K_{\infty}^{l-1}/{\mathcal{O}'}_{\infty}^{l-1}$.
There exist $c_{x,n}$ and $d_{x,n}$ such that $c < c_{x,n} < d_{x,n} <d$ and
$\{\widetilde{F}_n(t) + {\mathcal{O}'}_{\infty}^{l-1}| c_{x,n} \leq t \leq d_{x,n}\}$ is of length $\frac{\varepsilon}{2}$ and contained in
an $\varepsilon$-neighborhood of $x$. In view of the definition of $\widetilde{F}_n$, there exists $\delta$ not depending on $x$ and $n$
such that $e^{2\pi n n_{l-1}}(d_{x,n} - c_{x,n}) \geq \delta$. Now, since $n_l > n_{l-1}$, choosing $n$ sufficiently large we get
that $e^{2\pi n n_{l}}(d_{x,n} - c_{x,n}) \geq A_\varepsilon$ completing the proof of the claim.\qed

\subsection{A refinement of Jacobson-Morozov lemma.}
We will need the following known lemma (see \cite[Lemma 3.1]{EinLind}):

\begin{lem}
\label{Tits}
Let $\mathbf{L}$ be a semisimple group over a field $F$ of characteristic $0$, $\mathbf{S}$ be a maximal $F$-split torus in $\mathbf{L}$,
$\alpha$ be an indivisible root with respect to $\mathbf{S}$ and $\mathbf{V}_{(\alpha)}$ be the corresponding to $\alpha$ root group. Denote
$$U = \big\{ \left(
\begin{array}{cc} 1&x\\ 0&1\\ \end{array}
\right): x \in F \big\} \ \textit{and} \ D = \big\{ \left(
\begin{array}{cc} y&0\\ 0&y^{-1}\\ \end{array}
\right): y \in F^* \big\}.$$
Let $a = \exp (\nu)$ where $\nu \in \mathfrak{g}_\alpha$ if $2\alpha$ is not a root or $\nu \in \mathfrak{g}_\alpha \cup \mathfrak{g}_{2\alpha}$
otherwise.
Then there exists an $F$-morphism $f: \mathbf{SL}_2 \rightarrow \mathbf{L}$ such that
$a \in f(U)$ and $f(D) \subset \mathbf{S}(F)$.
\end{lem}

\subsection{Actions of epimorphic subgroups on homogeneous spaces in $\mathcal{S}$-adic setting.} Recall that $\G$ is a $K$-isotropic
semisimple $K$-group, $\mathcal{S} \supset \mathcal{S}_\infty$ and $G = G_\infty \times G_f$ where $G_\infty = \prod_{v \in \mathcal{S}_\infty}G_v$
and $G_f = \prod_{v \in \mathcal{S}_f}G_v$. Let $H$ be a closed subgroup of $G_\infty$ which have finite index in its Zariski closure. Recall that a subgroup
$B$ of $H$ is called \textit{epimorphic} if all $B$-fixed vectors are $H$-fixed for every rational linear representation of $H$. For example, the parabolic
subgroups in $H$ are epimorphic.
Also, if $f : K_v \rightarrow G_v, v \in \mathcal{S}$, is a $K_v$-rational homomorphism then $\{f(t): t \in K_1\}$ is called $1$-parameter unipotent
subgroup of $G_v$ (or $G$).

In the case when $\mathcal{S} = \mathcal{S}_\infty$ the following proposition is proved in \cite[Theorem 1]{Shah-Weiss}.

\begin{prop}\label{epimorphic subgroups} Let $H$ be a subgroup of $G_\infty$ generated by $1$-parameter unipotent subgroups and
$B$ be an epimorphic subgroup of $H$. Then any closed $B$-invariant subset of $G/\Gamma$ is $H$-invariant.
\end{prop}

{\bf Proof.} It is enough to prove that $\overline{B\pi(g)}$ is $H$-invariant for every $g \in G$. Since $g^{-1}Bg$ is an epimorphic subgroup of $g^{-1}Hg$,
the proof is easily reduced to the case when $g = e$.
Let $G_{f,n}$ be a decreasing sequence of
open compact subgroups of $G_f$ such that $\bigcap_n G_{f,n} = \{e\}$. Let $G_n = G_\infty \times G_{f,n}$ and $\Gamma_n = \Gamma \cap G_n$. Let
$\phi_n: G_n \rightarrow G_\infty$ be the natural projection and $\Gamma_{n,\infty} = \phi_n(\Gamma_n)$. Since $\Gamma_n$ is a lattice in $G_n$
 and $G_{f,n}$ is compact, $\Gamma_{n,\infty}$ is a lattice in $G_\infty$.
It follows from \cite[Theorem 1]{Shah-Weiss} that
$$
\overline{B\Gamma_{n,\infty}} = \overline{H\Gamma_{n,\infty}}
$$
for every $n$. Since $H$ is generated by $1$-parameter unipotent subgroups, in view of \cite{Ratner}, there exists a connected subgroup
$L_n$ of $G_\infty$ which contains $H$,
$$
\overline{H\Gamma_{n,\infty}} = {L_n\Gamma_{n,\infty}},
$$
and $L_n \cap \Gamma_{n,\infty}$ is a lattice in $L_n$. Since $G_{n+1}$ has finite index in $G_n$, $\Gamma_{n+1,\infty}$ has finite index in $\Gamma_{n,\infty}$.
Now using the connectedness of $L_n$ we get that all $L_n$ coincide, i.e. $L_n = L$.
So, $\phi_n(\overline{B\Gamma_n}) = \overline{B\Gamma_{n,\infty}} = L\Gamma_{n,\infty}$. Therefore for every
$x \in L$ there exists $a_n \in G_{f,n}$ such that $xa_n \in \overline{B\Gamma_{n}}$. Since $\{xa_n\}$
converges to $x$ in $G$, we get that $x \in \overline{B\Gamma}$, i.e. $L \subset \overline{B\Gamma}$.
Hence, in view of the inclusions $B \subset H \subset L$,
$$
\overline{B\Gamma} = \overline{H\Gamma},
$$
which proves our contention.
\qed

\subsection{Proof of Theorem \ref{r>2}.} We keep the assumptions from section 6.3 and suppose that $K_1 = \R$ or $\C$.
By Propositions \ref{prop4} and \ref{prop5} there exist $h \in \G(K)$, $\omega \in \mathcal{N}_G(T)$, and a nontrivial
defined over $K$ unipotent subgroup $\mathbf{U}$ of $\G$ such that if $U = \mathbf{U}(K_\mathcal{S})$ and $U_1 = \mathbf{U}(K_1)$ then
$\omega U \pi(h) \subset \overline{T\pi(g)}$ and $\mathcal{Z}_{T_1}(U_1)$ is finite. Shifting the orbit $T\pi(g)$ from the left by $\omega^{-1}$
we reduce the proof of the theorem to the case when $\omega = e$. Let $P_1$ be the maximal subgroup of $G_1$ with the properties:
$P_1$ is generated by $1$-parameter unipotent
subgroups of $G_1$ and ${(P_1 \times \{e\}\times \cdots \times \{e\})\pi(h)} \subset \overline{T\pi(g)}$.
Note that $U_1 \subset P_1$ and, therefore, $\mathcal{Z}_{T_1}(P_1)$ is finite. Since the projection
of $\te(\mathcal{O})$ into $T_1$ is Zariski dense and the stabilizer of $\pi(h)$ in $\te(\mathcal{O})$ has finite index in $\te(\mathcal{O})$,
$P_1$ is normalized by $T_1$. Put $X = \overline{T(P_1 \times \{e\}\times \cdots \times \{e\})\pi(h)}$. It is clear that $X \subset \overline{T\pi(g)}$
and that every $1$-parameter unipotent subgroup of $G_1$ which fixes $X$ (after its natural embedding in $G$) is contained in $P_1$.

Let us prove that $P_1$ is semisimple.
Suppose on the contrary that the unipotent radical $\mathcal{R}_u(P_1)$ of $P_1$ is not trivial.
Hence, there exists $a \in \mathcal{R}_u(P_1), a \neq e,$ such that $a = \exp (\nu)$ where $\nu \in \mathfrak{g}_\alpha$ for some root $\alpha$ of
$G_1$ with respect to $T_1$. By Lemma \ref{Tits} there exists a $K_1$-morphism $f: \mathrm{{SL}}_2(K_1) \rightarrow G_1$ such that
$a \in f(V)$ and $f(D) \subset T_1$ where $V$ is the subgroup of the upper triangular unipotent matrices in $\mathrm{{SL}}_2(K_1)$ and $D$ is the
subgroup of the diagonal matrices in $\mathrm{{SL}}_2(K_1)$. Denote by $B$ the subgroup of $P_1$ generated by $f(V)$ and $f(D)$. Then $B$ is
an epimorphic subgroup of $f(\mathrm{{SL}}_2(K_1))$ which fixes $X$. It follows from Proposition \ref{epimorphic subgroups} that $X$ is fixed
by $f(\mathrm{{SL}}_2(K_1))$ too. Therefore $f(\mathrm{{SL}}_2(K_1)) \subset P_1$.
Since $f(U) \subset \mathcal{R}_u(H_1)$ and the unipotent subgroups of $f(\mathrm{{SL}}_2(K_1))$ are conjugated
we get that $f(\mathrm{{SL}}_2(K_1)) \subset \mathcal{R}_u(H_1)$ which is a contradiction because $f(\mathrm{{SL}}_2(K_1))$
is not solvable.

Let $\widetilde{P}_1$ be the Zariski closure of $P_1$ in $G_1$. Then $P_1$ has finite index in $\widetilde{P}_1$ and $T_1 \subset P_1$ because
$T_1$ normalizes $\widetilde{P}_1$, $\widetilde{P}_1$ is semisimple  and
$\mathcal{Z}_{T_1}(\widetilde{P}_1)$ is finite. By \cite[Theorems 1 and 3]{Toma4} there exists unique connected $K$-subgroup $\mathbf{H}$ of $\G$
with the following properties: $\overline{{(P_1 \times \{e\}\times \cdots \times \{e\})\pi(h)}} = H\pi(h)$ where $H$ is a subgroup of
finite index in $\mathbf{H}(K_{\mathcal{S}})$, $P_1$ is contained in $H_1 = \mathbf{H}(K_1)$, and for each proper normal $K$-algebraic subgroup $\mathbf{Q}$
of $\mathbf{H}$ there exists $1 \leq i \leq r$ such that $(\mathbf{H}/\mathbf{P})(K_i)$ contains a unipotent element different from the identity.
Since $\mathcal{R}_u(H_1) \subset P_1$ and $P_1$ is semisimple we get that $\mathbf{H}$ is semisimple. In particular, $H\pi(h)$ is closed.
Also, a subgroup of finite index in
$\te(\mathcal{O})$ fixes $\overline{{(P_1 \times \{e\}\times \cdots \times \{e\})\pi(h)}}$ which implies that $\te$ normalizes $\mathbf{H}$. Therefore
$\te \subset \mathbf{H}$ and we may assume that $T \subset H$. So, $\overline{T\pi(g)}$ contains the closed $T$-invariant orbit $H\pi(h)$.
\qed

\subsection{Proof of Theorem \ref{cor5}.} We suppose that $\G = \mathbf{SL}_{n}, n \geq 2$.
As usual, $\mathbf{SL}_{n} = \mathbf{SL}(\mathbf{W})$ where $\mathbf{W}$ is the $K$-vector space
with $\mathbf{W}(K) = K^{n}$ and $\mathbf{W}(\OO) = \OO^{n}$.

Theorem \ref{cor5} follows from the next proposition.

\begin{prop}\label{borel- de siebenthal} Let $\mathbf{H}$ be a Zariski connected reductive $K$-subgroup of $\mathbf{SL}_{n}$
containing the subgroup of diagonal matrices $\te$ of $\mathbf{SL}_{n}$ and let $\mathbf{W} = \mathbf{W}_1 \oplus \cdots \oplus \mathbf{W}_l$ where
$\mathbf{W}_i$ are irreducible $\mathbf{H}$-subspaces defined over $K$. Then
\begin{enumerate}
\item[(a)] $\mathbf{H} = \{g \in \mathbf{SL}_{n}: g\mathbf{W}_i = \mathbf{W}_i \ \mathrm{for} \ \mathrm{all} \ i\}$ and
$\mathcal{D}(\mathbf{H}) = \mathbf{SL}(\mathbf{W}_1) \times \cdots \times \mathbf{SL}(\mathbf{W}_l)$;
\item[(b)] if $\mathbf{H}'$ is a reductive $K$-subgroup of $\mathbf{SL}_{n}$ such that $\mathbf{H}' \supset \mathbf{H}$ and
$\mathrm{s.s.rank}_K(\mathbf{H})=\mathrm{s.s.rank}_K(\mathbf{H}')$ then $\mathbf{H} = \mathbf{H}'$.
\end{enumerate}
\end{prop}

{\bf Proof.} The proof is based on the following observation. Let $\mathbf{L} = \mathbf{L}_1 \times \cdots \times \mathbf{L}_s$ be a direct
product of simple algebraic $K$-groups and let $\rho: \mathbf{L} \rightarrow \mathbf{SL}(\mathbf{V})$ be an irreducible $K$-representation
with finite kernel.
It is well-known that $\rho$ is a tensor product of irreducible representations of $\mathbf{L}_i, 1 \leq i \leq s$. In view of the description
in \cite[Table 2]{Bourbaki} of the dimensions of the irreducible representations of the simple algebraic groups, we have that
$\dim (\mathbf{V}) \geq \mathrm{rank}_K (\mathbf{L}) + 1$ and $\dim (\mathbf{V}) = \mathrm{rank}_K (\mathbf{L}) + 1$ if and only if $\rho$ is a
$K$-isomorphism.

The group $\mathbf{H}$ as in the formulation of the proposition is an almost direct product of its center $\mathbf{Z}$ and $\mathcal{D}(\mathbf{H})$.
Hence, every $\mathbf{W}_i$ is an irreducible
$\mathcal{D}(\mathbf{H})$-subspace, $\dim \mathbf{Z} \leq l-1$ and,
since $\mathcal{D}(\mathbf{H}) \subset \mathbf{SL}(\mathbf{W}_1) \times \cdots \times \mathbf{SL}(\mathbf{W}_l)$,
$\mathrm{rank}_K (\mathcal{D}(\mathbf{H})) \leq \overset{l}{\underset{i = 1}{\sum}}(\dim (\mathbf{W}_l) - 1) = n -l$. But
$\mathrm{rank}_K (\mathbf{H}) = \mathrm{rank}_K (\mathcal{D}(\mathbf{H})) + \dim \mathbf{Z} = n-1$. Therefore, $\dim \mathbf{Z} = l-1$ and
$\mathrm{s.s.rank}_K (\mathcal{D}(\mathbf{H})) = \overset{l}{\underset{i = 1}{\sum}}(\dim (\mathbf{W}_i) - 1)$.
It follows from the above observation that
$\mathcal{D}(\mathbf{H}) = \mathbf{SL}(\mathbf{W}_1) \times \cdots \times \mathbf{SL}(\mathbf{W}_l)$ and
$\mathbf{H} = \{g \in \mathbf{SL}_{n}: g\mathbf{W}_i = \mathbf{W}_i \ \mathrm{for} \ \mathrm{all} \ i\}$, proving (a). Note that every
irreducible $\mathcal{D}(\mathbf{H}')$-subspace is $\mathcal{D}(\mathbf{H})$-invariant and the center of $\mathbf{H}'$ is contained
in $\mathbf{Z}$. This implies (b).\qed

\medskip

Theorem \ref{cor5} is deduced from Theorem \ref{r>2} and Proposition \ref{borel- de siebenthal} as follows. Let $\mathbf{H}_1$, $\mathbf{H}_2$ and $\te$ be as in the
formulation of Theorem \ref{r>2}.
It follows from Proposition \ref{borel- de siebenthal} that $\mathbf{H}_1 = \mathbf{H}_2 = \mathbf{H}$.
Let $H$ be a subgroup of finite index in $\mathbf{H}(K_\mathcal{S})$ containing $\te(K_\mathcal{S})$. Since $\mathrm{SL}_t(F)$ does
not contain subgroups of finite index whenever $F$ is a field of
characteristic $0$, we get that
$H = \mathbf{H}(K_\mathcal{S})$ which implies Theorem \ref{cor5}.

\medskip

\section{A number theoretical application}
\label{spit forms}

\subsection{Reduction of the proof of Theorem \ref{application} to the case $m = n$.} The reduction is based on the following general
\begin{prop} \label{reduction}
 Let $M_i$, where $1 \leq i \leq r$ and $r > 1$, be subsets of the vector space $K^n$ each of them consisting of $m$ linearly independent
 vectors. Suppose that there exist $i \neq j$ and $\vec{w} \in M_i$ such that $K\vec{w} \cap M_j = \emptyset$.
 Then there exists a linear map
 $\phi: K^n \rightarrow K^m$ such that every $\phi(M_i)$ consists of $m$ linearly independent vectors and
 $K\phi(\vec{w}) \cap \phi(M_j) = \emptyset$.
\end{prop}

{\bf Proof.} Let $\mathcal{M}$ be the subset of $\mathrm{End}_K(K^n,K^m)$ consisting of all $\phi$ as in the formulation
of the proposition. One can prove by a standard argument that $\mathcal{M}$ is a Zariski open non-empty subset of $\mathrm{End}_K(K^n,K^m)$
which proves the proposition. \qed

With the notation from the formulation of Theorem \ref{application}, let $m < n$. We identify $K^n$ with $K x_1 + \cdots + K x_n$ and put
$M_v = \{l_{1}^{(v)}(\vec{x}), \ldots,
l_{m}^{(v)}(\vec{x}) \}$, $v \in \mathcal{S}$. By the assumptions of the theorem there exist $v_1$ and $v_2 \in \mathcal{S}$ and $1 \leq i \leq m$
such that
$Kl_{i}^{(v_1)}(\vec{x}) \cap M_{v_2} = \emptyset$. There exists a basis $\vec{y} = (y_1,\ldots,y_n)$ of
$K x_1 + \cdots + K x_n$ such that the map
$\phi$ as in the formulation of Proposition \ref{reduction} is given by $\phi(a_1y_1 + \cdots + a_n y_n) = b_1y_1 + \cdots + b_m y_m$
where $b_s$ depend linearly on $a_t$. Every $l_{i}^{(v)}(\vec{x})$ is a linear form on $y_1,\ldots,y_n$ denoted by
$\lambda_{i}^{(v)}(\vec{y})$. Put $\widetilde{l}_{1}^{(v)}(y_1, \cdots, y_m) = \lambda_{i}^{(v)}(y_1, \cdots, y_m, 0, \cdots, 0)$ and
$\widetilde{f}_v(y_1, \cdots, y_m) = \overset{m}{\underset{i=1}{\prod}} \widetilde{l}_{i}^{(v)}(y_1, \cdots, y_m)$, $v \in \mathcal{S}$.
Then $\widetilde{l}_{1}^{(v)}, \cdots, \widetilde{l}_{m}^{(v)}$
are linearly independent over $K$ and $\widetilde{f}_{v_1}$ is not proportional to $\widetilde{f}_{v_2}$.
So, the validity of the theorem for ${f} \in K_{\mathcal{S}}[x_1, \cdots, x_n]$
follows from its validity for $\widetilde{f} = (\widetilde{f}_v)_{v \in \mathcal{S}} \in K_{\mathcal{S}}[y_1, \cdots, y_m]$. \qed

In the framework of Theorem \ref{application}, it is a natural problem to understand  the distribution of $f(\OO^n)$ in $K_{\mathcal{S}}$.
Presumably, it is a matter of uniform distribution.

\subsection{Proof of Theorem \ref{application}.}
Let $\G = \mathbf{SL}_{n}$, $G = \mathbf{SL}_{n}(K_{\mathcal{S}})$ and
$\Gamma = \mathbf{SL}_{n}(\OO)$. The group $G$ is acting on $K_{\mathcal{S}}[
\vec{x}]$ according to the law $(\sigma \phi)(\vec{x}) =
\phi(\sigma^{-1}\vec{x})$, where $\sigma \in G$ and $\phi \in K_{\mathcal{S}}[
\vec{x}]$. We denote $f_0(\vec{x}) = x_1x_2...x_{n}
$. Let $f(\vec{x})$
be as in the formulation of the theorem with $m = n$. There exists $g = (g_v)_{v \in \mathcal{S}} \in G$ such that
every $g_v \in \G(K)$ and $f(\vec{x}) = \alpha(g^{-1} f_0)(\vec{x})$ where $\alpha \in K_{\mathcal{S}}$. Since
$f_v(\vec{x}), v_n \in S,$ are not pairwise proportional the orbit $T\pi(g)$ is locally divergent but non-closed
(Theorem \ref{ldo}). Note that $f(\vec{x}) = \alpha(w g f_0)(\vec{x})$ for every $w \in \mathcal{N}_G(T)$. In view of
Theorem \ref{cor5} there exist a reductive $K$-subgroup $\mathbf{H}$ with $\te \varsubsetneq \mathbf{H}$ and $\sigma \in \G(K)$
such that  $\overline{T\pi(g)} = H\pi(\sigma)$ where $H = \mathbf{H}(K_{\mathcal{S}})$. By  Proposition \ref{borel- de siebenthal} there
exists a direct sum decomposition
$\mathbf{W} = \mathbf{W}_1 \oplus \cdots \oplus \mathbf{W}_l$ such that
$\mathbf{H} = \{g \in \mathbf{SL}_{n}: g\mathbf{W}_i = \mathbf{W}_i \ \mathrm{for} \ \mathrm{all} \ i\}$ and at least one of the subspaces
$\mathbf{W}_i$, say $\mathbf{W}_{i_\circ}$,  has dimension $> 1$.
(Recall that $\mathbf{SL}_{n} = \mathbf{SL}(\mathbf{W})$.) 
Pick $a = (a_v)_{v \in \mathcal{S}} \in K_{\mathcal{S}}$ with $a_v \neq 0$ for all $v$. Using $\dim \mathbf{W}_{i_\circ} > 1$,
we find $\vec{z} \in \OO^{n}$ and $h = (h_v)_{v \in \mathcal{S}} \in H$ such that
$f_0(h\sigma(\vec{z})) = a$.
Since $H\sigma \Gamma = \overline{Tg\Gamma}$ there exist $t_i \in T$ and $\gamma_i \in \Gamma$ with
$$
\lim_i t_ig\gamma_i = h\sigma.
$$
Therefore
$$
\lim_i f(\gamma_i \vec{z}) = a,
$$
proving the theorem. \qed

\section{Examples}
In this section we provide examples showing that $\overline{T\pi(g)}$ in the formulation of Theorem \ref{r>2} might be not homogeneous and that
the claims of Theorem \ref{cor5}
and Theorem \ref{application} are not true for a $\mathrm{CM}$-field $K$. For simplicity, we will suppose that $\mathcal{S} = \mathcal{S}_{\infty}$.
So, let $K$ be a $\mathrm{CM}$-field, that is, $K = F(\sqrt{-d})$, where $F$ is a totally real number field, $d \in F$ and $d > 0$ in every archimedean
completion of $F$. We denote in the same way the archimedean places of $F$ and their (unique) extensions to $K$. So, $K_v = \C$ and $F_v = \R$
for all $v \in \mathcal{S}$.
Also, let $\OO_F$ (resp. $\OO_K$) be the ring of integers of $F$ (resp. $K$). Recall that
$\OO_F$ (resp. $\OO_K$) is a lattice in $F_\mathcal{S} = {\underset{v \in \mathcal{S}}{\prod}}F_v$ (resp. $K_\mathcal{S} = {\underset{v \in \mathcal{S}}{\prod}}K_v$).

\subsection{Restriction of scalars functor for $\mathrm{CM}$-fields.}

Denote by $\mathbf{G}$ the group $\mathbf{SL}_{2}$ considered as a $K$-algebraic group. Let $^-: K \rightarrow K$ be the non-trivial automorphism of $K/F$. For every $v \in \mathcal{S}$ we keep the same notation $^-$ for the complex conjugation of $K_v = \C$ and for the group automorphism
 $\SL_2(K_v) \rightarrow \SL_2(K_v), \left(
\begin{array}{cc} x&y\\ z&t\\ \end{array} \right) \mapsto \left(
\begin{array}{cc} \overline{x}&\overline{y}\\ \overline{z}&\overline{t}\\ \end{array} \right)$. There exists a simple $F$-algebraic group of $F$-rank $1$,
denoted by $R_{K/F}(\mathbf{G})$, and a $K$-morphism $p: R_{K/F}(\mathbf{G}) \rightarrow \G$ such that the map
$(p, \overline{p}): R_{K/F}(\mathbf{G}) \rightarrow
\G \times \G, g \mapsto (p(g), \overline{p(g)}),$ is a $K$-isomorphism of $K$-algebraic groups and $p(R_{K/F}(\mathbf{G})(F)) = \G(K)$. 
The pair $(R_{K/F}(\mathbf{G}), p)$ is uniquely defined by the above
properties up to an $F$-isomorphism and the $F$-algebraic group $R_{K/F}(\mathbf{G})$ is obtained from the $K$-algebraic group $\G$ via
 the restriction scalars functor $R_{K/F}$. (We refer to \cite[6.17-6.21]{Borel-Tits} or \cite[1.3]{W2} for the general definition and basic properties
of $R_{K/F}$.)

Given $v \in \mathcal{S}$, the isomorphism $R_{K/F}(\mathbf{G})(F) \rightarrow \G(K), g \mapsto p(g)$ admits a unique extension to an isomorphism $R_{K/F}(\mathbf{G})(F_v) \rightarrow G_v$ denoted by $p_v$. Let $p_{\mathcal{S}}$ be the direct product of all $p_v, v \in \mathcal{S}$. Further on $R_{K/F}(\mathbf{G})(F_{\mathcal{S}})$ will be identified with $G$ via the isomorphism $p_{\mathcal{S}}$.
Let
$\T$ be the subgroup of the diagonal matrices in $\G$. Under the above identification $\Gamma = \G(\OO_K) = R_{K/F}(\mathbf{G})(\OO_F)$
and $T = \T(K_{\mathcal{S}}) = R_{K/F}(\mathbf{T})(F_{\mathcal{S}})$. For every $v \in \mathcal{S}$ we have that $T_v = \T(K_{v})$ is the the group of
complex diagonal matrices in $G_v (= \SL_2(\C))$. The $F$-torus $R_{K/F}(\T)$ is not split and contains a maximal $1$-dimensional $F$-split torus $\T_F$.
Note that $\T_F(F_{v})$, $v \in \mathcal{S}$,  is the the group of
\textit{real} diagonal matrices in $T_v$. Denote $T_\R = \T_F(F_{\mathcal{S}})$.  Then $T = T_{\R} \cdot N$ where $N$ is a compact group.

\subsection{Non-homogeneous $T$-orbits closures when $r > 2$}

We continue to use the notation and the assumptions from \S8.1.
Also, let $u^-(x) = \left(
\begin{array}{cc} 1&0\\ x&1\\ \end{array} \right)$ and $u^+(x) = \left(
\begin{array}{cc} 1&x\\ 0&1\\ \end{array} \right)$, $\mathcal{S} = \{v_1, \cdots, v_r\}$ and $G = G_{v_1} \times \cdots \times G_{v_r}$.

\begin{thm}\label{contre-example} Suppose that $r > 2$. Let
$g = (u^-(\beta)u^+(\alpha), e, \cdots, e)
\in G$ where $\alpha \in F^*$ and $\beta \in K \setminus F$. Then the following holds:
\begin{enumerate}
\item[(a)] Each of the orbits ${T\pi(g)}$ and ${T_\R \pi(g)}$ is not dense in $G/\Gamma$,
\item[(b)] Each of the sets $\overline{T\pi(g)} \setminus T\pi(g)$ and $\overline{T_\R\pi(g)} \setminus T_\R\pi(g)$ is not contain in a union of countably
many closed orbits of proper subgroups of $G$.
\end{enumerate}
In particular, each of the closures $\overline{T\pi(g)}$ and $\overline{T_\R \pi(g)}$ is not homogeneous.
\end{thm}

{\bf Proof.}
A direct calculation shows that $u^-(\beta)u^+(\alpha) = du^+(\alpha_1)u^-(\beta_1)$ where $\alpha_1 = (1 + \alpha\beta)\alpha$, $\beta_1 = (1+ \alpha \beta)^{-1} \beta$ and $d = \left(
\begin{array}{cc} (1+\alpha\beta)^{-1}&0\\ 0 &1+\alpha\beta\\ \end{array}
\right)$. Since $\beta \in K \setminus F$ we get that $\beta_1 \in K \setminus F$.

Define subgroups $L_1$ and $L_2$ of $G$ as follows. Put $L_1 = \mathbf{SL}_2(F_{\mathcal{S}})$ and
$L_2 = \{\left(
\begin{array}{cc} x&y\beta_1^{-1}\\ z\beta_1&t\\ \end{array}
\right) \in G: x,y,z,t \in F_{\mathcal{S}}\}$. The group $L_1 \cap \Gamma$ is commensurable with
$\SL_2(\OO_F)$ and, therefore, is a lattice in $L_1$. Hence $L_1\pi(e)$ is closed.
Since the map $\G \rightarrow \G, \left(
\begin{array}{cc} x&y\\ z&t\\ \end{array}
\right) \mapsto \left(
\begin{array}{cc} x&y\beta_1^{-1}\\ z\beta_1&t\\ \end{array}
\right)$, is a $K$-isomorphism we get that $L_2\pi(e)$ is closed too.

It follows from the definitions of $L_1$ and $L_2$ that
\begin{equation}
\label{subset}
{T_\R \pi(g)} \subset \underset{0\leq
\mu\leq1}{\bigcup}\{(u^-(\mu\beta), \cdots, e)L_1\pi(e)\} \bigcup
\underset{0\leq \nu\leq1}{\bigcup}\{(d
\cdot u^+(\nu\alpha_1), \cdots, e)L_2\pi(e)\}.
\end{equation}
Since the right hand side of (\ref{subset}) is a proper closed subset of $G/\Gamma$, $\overline{T_\R \pi(g)} \neq G/\Gamma$.
Also, it is easy to see that the shift of the right hand side of (\ref{subset}) by the compact group $N$ (defined at the end of \S8.1) remains
 a proper subset of $G/\Gamma$.
Since $T = T_{\R} \cdot N$, $\overline{T \pi(g)} \neq G/\Gamma$, completing the proof of (a).

Let ${U_1}^+$ be the subgroup of all upper triangular unipotent matrices in ${L}_1$ and ${U_2}^-$ be the group of all lower triangular unipotent matrices
in ${L}_2$.
It follows from Proposition \ref{prop5} and Proposition \ref{integers1}(2a)
that $\overline{T_\R^\bullet \pi(g)} \supset U_1^+\pi(e) \cup (d,e, \cdots, e)U_2^-\pi(e)$. In view of Proposition \ref{epimorphic subgroups}
we have
\begin{equation}
\label{subset1}
\overline{T_\R^\bullet \pi(g)} \supset L_1\pi(e) \cup (d,e, \cdots, e)L_2\pi(e).
\end{equation}

Fix a real transcendental number $a$. Using again Proposition \ref{integers1} we get that $\overline{T_\R^\bullet\pi(g)}$ contains $\pi(\widetilde{g})$
where $\widetilde{g} = (u^-(a^{-2}\beta)u^+(a^{2}\alpha), \cdots, e)$. Note that
\begin{equation}
\label{subset1+}
\overline{T_\R^\bullet  \pi(\widetilde{g})} = \overline{T_\R^\bullet  \pi({g})}.
\end{equation}
Suppose that $\pi(\widetilde{g}) \in T\pi(g)$. Then there exist $t = \left(\begin{array}{cc} \tau&0\\ 0&\tau^{-1}\\ \end{array}
\right) \in \SL_2(\C)$ and $m \in \SL_2(K)$ such that $u^-(\beta)u^+(\alpha)m = t u^-(a^{-2}\beta)u^+(a^{2}\alpha)$.
The upper left coefficient of $t u^-(a^{-2}\beta)u^+(a^{2}\alpha)$ is equal to $\tau$ and $u^-(\beta)u^+(\alpha)m \in \SL_2(K)$. Hence $\tau \in K$. On the other hand, the upper right coefficient of $t u^-(a^{-2}\beta)u^+(a^{2}\alpha)$ is equal to $\tau a^2 \alpha$. Therefore
$a$ is an algebraic number which is a contradiction. We have proved that
$$
\pi(\widetilde{g}) \in \overline{T_\R^\bullet \pi(g)} \setminus T\pi(g).
$$
Since
$$
\overline{T_\R^\bullet \pi(g)} \setminus T\pi(g) \subset (\overline{T_\R^\bullet \pi(g)} \setminus T_\R^\bullet  \pi(g)) \cap (\overline{T\pi(g)} \setminus T\pi(g)),
$$
in order to prove (b) it is enough to show that if
$\overline{T_\R^\bullet  \pi({g})} \setminus T\pi(g) \subset \overset{\infty}{\underset{i=1}{\bigcup}} Q_i\pi(h_i)$,
where $Q_i$ are connected closed subgroups of $G$ and $Q_i\pi(h_i)$ are closed orbits, then one of the subgroups $Q_i$ is equal to $G$.
It follows from Baire's category theorem, applied to a compact neighborhood of $\pi(\widetilde{g})$ in $\overline{T_\R^\bullet\pi(\widetilde{g})}$,
that there exists $i_\circ$ such that $T_\R^\bullet  \subset Q_{i_\circ}$ and
 $\pi(\widetilde{g}) \in Q_{i_\circ}\pi(h_{i_\circ})$.
 Using (\ref{subset1}) and (\ref{subset1+}), we obtain that
 $L_1 \cup (d,e, \cdots, e)L_2(d,e, \cdots, e)^{-1} \subset Q_{i_\circ}$. Since $\beta_1 \in K \setminus F$ it follows from the definitions of $L_1$
 and $L_2$ that $Q_{i_\circ}$ contains $\{e\} \times \cdots \times \G(K_{v_r})$. 
But $\Gamma$ is an irreducible lattice in $G$. Therefore $(\{e\} \times \cdots \times  \G(K_{v_r}))\cdot \Gamma$ is dense in $G$ and
 $Q_{i_\circ} = G$.
\qed

\medskip

\textsf{Remark.}
The orbit $T\pi(g) \subset \mathbf{SL}_{2}(K_\mathcal{S})/\mathbf{SL}_{2}(\OO_K)$
provides an example showing that Theorem \ref{cor5} is not valid for $\mathrm{CM}$-fields.
On the other hand, $T_\R\pi(g) \subset R_{K/F}(\mathbf{\mathbf{SL}_{2}})(F_\mathcal{S})/R_{K/F}(\mathbf{\mathbf{SL}_{2}})(\OO_F)$
 provides an example showing that $\overline{T\pi(g)}$ as in the formulation of Theorem \ref{r>2} is not always homogeneous.

\subsection{Values of decomposable forms when $\#\mathcal{S} = 2$ or $\#\mathcal{S} \geq 2$ and $K$ is a $\mathrm{CM}$-field.}
Let us provide the necessary counter-examples showing that
the assertion of Theorem \ref{application} does not hold if $\#\mathcal{S} = 2$ or $K$ is a $\mathrm{CM}$-field and $\#\mathcal{S} \geq 2$.

We keep the notation $f(\vec{x})$, $f_v(\vec{x})$ and $l_i^{(v)}(\vec{x})$ as in the formulation of Theorem \ref{application}.
We will assume that $m = n = 2$.

The following is a particular case of \cite[Theorem 1.10]{Toma3}:

\begin{thm}\label{two valuations} Let $\#\mathcal{S} = 2$. Then $\overline{f(\OO_K^2)} \cap K_{\mathcal{S}}^*$ is a countable set. In particular, $\overline{f(\OO_K^2)}$ is not dense in $K_{\mathcal{S}}^*$.
\end{thm}

Remark that in the formulation of Theorem \ref{two valuations} $K$ is not necessarily a $\mathrm{CM}$-field.

\begin{thm}\label{example CM field} Let $K$ be a $\mathrm{CM}$-field which is a quadratic extension of
a totally real field $F$. We suppose that $\#\mathcal{S} \geq 2$ and that all $l_i^{(v)}$ are with coefficients from $F$.
Then there exists a real $C > 0$ such that for every $\vec{z} \in \OO_K^2$ either
\begin{equation}
\label{Artin1}
\underset{v \in \mathcal{S}}{\prod}|f_{v}(\vec{z})|_v \in C\cdot \N,
\end{equation}
or there exists $\alpha \in \C^*$ such that
\begin{equation}
\label{Artin2}
 f_v(\vec{z}) \in \alpha \cdot\R  \ \mathrm{for} \ \mathrm{all} \ v \in \mathcal{S}.
\end{equation}
In particular, $\overline{f(\OO^2)}$ is not dense in $K_{\mathcal{S}}$.
\end{thm}

{\bf Proof.} Choose $d \in \OO_F$ such that $K = F(\sqrt{-d})$.
Let $\sigma_i$, $1 \leq i \leq r$, be the set of all nontrivial morphisms of the field $F$ into $\C$.
Given $\sigma_i$, we will keep the same notation for its extension to the morphism from $K$ into $\C$ which maps
 $\sqrt{-d}$ to $\sqrt{-\sigma_i(d)}$. The $\textit{normalized}$ archimedean valuation of $K$ corresponding to $\sigma_i$ is defined by
 $| x |_{v_i}:= \| \sigma_i(x) \|^2$, where $x \in K$ and $\| \cdot \|$ is the usual norm on $\C$.
Also, recall that
$\mathrm{N}_{K/\Q}(x) := \underset{i}{\prod}| x |_{v_i}$ is the algebraic norm of $x$ and $\mathrm{N}_{K/\Q}(x) \in \N$ if
$x \in \OO_K$ (see \cite[ch.2, Theorem 11.1]{Cassels-Fröhlich}).

Let $l_{i}^{(v_j)}(x_1,x_2) = h_{i1}^{(j)}x_1 + h_{i2}^{(j)}x_2$ where $j \in \{1, \cdots, r\}$ and $i \in \{1, 2\}$.
Put $h^{(j)} := \left(
\begin{array}{cc} h_{11}^{(j)}&h_{12}^{(j)}\\ h_{21}^{(j)}&h_{22}^{(j)}\\ \end{array} \right)$.
We have $f_{j} = l_{1}^{(v_j)}\cdot l_{2}^{(v_j)}$. Multiplying $f_{j}$ by appropriate elements from $F^*$
we suppose without loss of generality that all $h^{(j)} \in \SL_2(F)$. Further on, if $\vec{w}_1 = (w_{11}, w_{12}) \in \C^2$ and $\vec{w}_2 = (w_{21}, w_{22}) \in \C^2$ we denote by $\det (\vec{w}_1,\vec{w}_2)$ the determinant of $\left(
\begin{array}{cc} w_{11}&w_{12}\\ w_{21}&w_{22}\\ \end{array} \right)$. Also, given $w = \left(
\begin{array}{cc} w_{11}&w_{12}\\ w_{21}&w_{22}\\ \end{array} \right)$ and
$\vec{a} = (a_1, a_2) \in \C^2$ we write $w(\vec{a}) = (w_{11}a_1 + w_{12}a_2, w_{21}a_1 + w_{22}a_2)$.

Let $\vec{z} = \vec{\gamma} + \sqrt{-d}\vec{\delta} \in \OO_K^2$ where $\vec{\gamma} = (\gamma_1, \gamma_2) \in F^2$ and
$\vec{\delta} = (\delta_1, \delta_2) \in F^2$. Denote by $l$ the index of  $\OO_F[\sqrt{-d}]$ in $\OO_K$. Then $(l\gamma_1, l\gamma_2) \in \OO_F^2$ and
$(l\delta_1, l\delta_2) \in \OO_F^2$. Since $\det(h^{(j)}(\vec{\gamma}), h^{(j)}(\vec{\delta})) = \det(\vec{\gamma}, \vec{\delta})$, we get
$$
\underset{j}{\prod}|\det(h^{(j)}(\vec{\gamma}), h^{(j)}(\vec{\delta}))|_{v_j} =
\underset{j}{\prod}|\det(\vec{\gamma}, \vec{\delta})|_{v_j} \in \frac{1}{l^{4r}}\N.
$$
We have
$$
\sigma_j(l_{i}^{(v_j)}(\vec{z})) = \sigma_j(l_{i}^{(v_j)}(\vec{\gamma}))+ \mathfrak{i}\sqrt{\sigma_j(d)}\sigma_j(l_{i}^{(v_j)}(\vec{\delta}))=
r_i^{(j)}e^{\mathfrak{i}\varphi_i^{(j)}}
$$
and
$$
|{f_{j}(\vec{z})}|_{v_j} = (r_1^{(j)}\cdot r_2^{(j)})^2,
$$
where $r_i^{(j)}$ (resp. $\varphi_i^{(j)}$) is the absolute value (resp. the argument) of the complex number $\sigma_j(l_{i}^{(v_j)}(\vec{z}))$.
A simple computation shows that
$$
|\det (h^{(j)}(\vec{\gamma}),h^{(j)}(\vec{\delta}))|_{v_j} = \frac{|f_{j}(\vec{z})|_{v_j}}{\sigma_j(d)}|\sin (\varphi_{1}^{(j)}-\varphi_{2}^{(j)})|^2.
$$
Therefore if $\vec{\gamma}$ and $\vec{\delta}$ are not proportional then
$$
\underset{j}{\prod}|f_{j}(\vec{z})|_{v_j} \geq \sqrt{\mathrm{N}_{K/\Q}(d)} \underset{j}{\prod}|\det(\vec{\gamma}, \vec{\delta})|_{v_j}  \in
\frac{\sqrt{\mathrm{N}_{K/\Q}(d)}}{l^{4r}}\N,
$$
proving (\ref{Artin1}). Let $\vec{\delta} = a \vec{\gamma}$, $a \in F$. Then $f_{v_j}(\vec{z}) =
(1+a \sqrt{-d})^2f_{v_j}(\vec{\gamma})$ for all $j$ which implies (\ref{Artin2}).\qed

\medskip

Acknowledgement. This work was partially supported by the IMI of BAS.

\end{document}